\documentclass[twocolumn]{autart}    % Enable this line and disable the 
\usepackage{graphicx}          % Include this line if your 
\newtheorem{proof}{Proof}
\newtheorem{example}{Example}

\usepackage{amsmath}
\usepackage{amssymb}
\usepackage{amsfonts}
\usepackage{mathtools}
\usepackage{color}
\usepackage{textcomp}

\usepackage[utf8]{inputenc}

\newcommand{\R}{\mathbb{R}}

\newcommand{\K}{\mathcal{K}}
\renewcommand{\P}{\mathcal{P}}
\renewcommand{\H}{\mathcal{H}}

\newcommand{\D}{\mathcal{D}^+}
\renewcommand{\S}{\mathcal{S}}
\newcommand{\<}{\leqslant}
\renewcommand{\>}{\geqslant}
\renewcommand{\(}{\left(}
\renewcommand{\)}{\right)}
\newcommand{\Hn}{{_H}}
\renewcommand{\vec}{\mathrm{vec}}
\newcommand{\vph}{\vphantom{\underline{\overline{\Omega}}}}
\newcommand{\vphS}{\vphantom{\bar\Omega}}
\newcommand{\Co}{\mathbf{C}\left(\R,\R^{n\times n}\right)}
\newcommand{\Cf}{\mathbf{C}^{(1)}\left([-H,0],\R^n\right)}
\newcommand{\PC}{\mathbf{PC}\left([-H,0],\R^n\right)}
\newcommand{\ph}{\varphi}
\newcommand{\eps}{\varepsilon}

\allowdisplaybreaks

\usepackage[round]{natbib}
\begin{document}

\begin{frontmatter}

\title{Lyapunov stability tests for linear time-delay systems\thanksref{footnoteinfo}} % Title, preferably not more 

\thanks[footnoteinfo]{This project was partially supported by Conacyt A1-S-24796.\\
Corresponding author: S. Mondié\\ smondie@ctrl.cinvestav.mx (S. Mondié)\\ ORCID(s): 0000-0001-7511-2910 (S. Mondié)}

\author[1]{Sabine Mondi\'e}%\ead{smondie@ctrl.cinvestav.mx},   
% \cormark[1]Add the 
\author[2]{Alexey Egorov},               % e-mail address 
\author[3]{Marco A. Gomez}  % (ead) as shown

\address[1]{Department of Automatic Control, CINVESTAV, Av. IPN 2508, Mexico City, 07360, Mexico} % Please supply                                              
\address[2]{Department of Control Theory, St. Petersburg State University, 7–9 Universitetskaya nab., St. Petersburg, 199034, Russia}             % full addresses
\address[3]{Department of Mechanical Engineering, DICIS, Universidad de Guanajuato, Salamanca, 36885, Guanajuato, Mexico}        % here.

                             % keyword list or with the 
                                          % help of the Automatica 
                                          % keyword wizard

\begin{abstract}                          % Abstract of not more than 200 words.
An overview of stability conditions in terms of the Lyapunov matrix for time-delay systems is presented. The main results and proof are presented in details for the case of systems with multiple delays. The state of the art, ongoing research and potential extensions to other classes of delay systems are discussed.
\end{abstract}

\end{frontmatter}

%==============================================================
%%%%%%%%%%%%%%%%%%%%%%%%%%%%%%%%%%%%%%%%%%%%%%%%%%%%%%%

%%%%%%%%%%%%%%%%%%%%%%%%%%%%%%%%%%%%%%%%%%%%%%%%%%%%%%%
\section{Introduction}

One of the most famous control theory results is the stability criterion (i.\,e., necessary and sufficient condition) in the Lyapunov framework. Indeed there is beauty in the fact that a simple matrix captures the complexity of the behavior of a linear dynamic system on an infinite time interval. This result stands out because of its simplicity, elegance, and power to address nonlinearities, parameter uncertainties, additive perturbations, and so on.

More precisely, the proof that the linear system
\begin{equation*}
	\dot{x}=Ax
\end{equation*}
is stable if and only if there exists a function of the form $v(x)=x^{T}Px$ that satisfies $v(x)>0$, $x\neq 0$, and $\frac{d}{dt}v(x)<0,$ $x\neq 0$, along the trajectories of the system established in the Lyapunov framework reduces, when the Lyapunov condition is satisfied, to solving for any positive definite matrix $W$, the algebraic equation
\begin{equation*}
	A^{T}P+PA=-W
\end{equation*}
for the unknown matrix $P$, and to verify its positive definiteness.

In the fifties, the foundations for the extension of the Lyapunov approach to time-delay systems were given by Razhumikhin and Krasovskii, with theorems linking stability to the existence of a function \citep{Razumikhin1956} or a functional \citep{Krasovskii1956}. The impact of the latter approach has been more significant, especially in the generation of sufficient stability conditions, consisting of the proposal of better functionals leading to sufficient stability conditions in the form of linear matrix inequalities. Such results dominate the literature, mainly due to the vast array of applications they allow. They range from the early proposals of \cite{Kolmanovskii1999,Niculescu2001} to refined conditions such as those of \cite{Fridman2014,SeuretGouaisbaut2015}, among many others.

The converse results of the theory of Krasovskii, which guarantee the existence of the functional when the system is stable, have been less popular, although the general form of the functional introduced by \cite{Repin1965} and \cite{Datko1972} has been a source of inspiration for the determination of sufficient stability conditions in \cite{Gu2001} and \cite{PEET2011}.

In the past decades, following the early works by \cite{InfanteandCastelan1978,Huang1989,Louisell2001}, a comprehensive coverage of the case of linear delay systems of retarded, distributed, and neutral types has been made by \cite{Kharitonov2013}. The core of the presented results includes the functional, which is obtained from the Cauchy formula and expressed in terms of the delay Lyapunov matrix, the construction of the Lyapunov matrix via the four properties called continuity, dynamic, symmetry and algebraic, and the proof of the existence of a positive quadratic lower bound whenever the system is stable.

These results are indeed a nontrivial extension of the delay-free case. They were successfully applied to robust stability analysis \citep{KharitonovandZhabko2003}, exponential estimates \citep{KharitonovHinrichsen2004}, computation of the $\H_2$ norm \citep{Jarlebringetal2011,SumachevaandKharitonov2014}, solution of the suboptimal control problem \citep{Santosetal2009} and design of predictor-based control schemes for state and input-delay systems \citep{Kharitonov2014}. Notably, these applications share the common assumption that the underlying system is exponentially stable. This raises the following question: \textit{is it possible to assess the system’s stability in this framework? Moreover, considering the analogy with the delay-free case, does a stability criterion expressed in terms of the delay Lyapunov matrix exist}?

An encouraging fact regarding the soundness of this aim was the finding of necessary and sufficient stability conditions in terms of the delay Lyapunov matrix for the scalar one delay equation. In this case, both the stability region and a simple analytical expression of the $1\times 1$ delay Lyapunov matrix are known. We were able to present in \cite{Mondie2012,EgorovMondie2013a} two sets of such conditions.

We have devoted the past years to answer this query for several classes of linear systems with delays. This contribution intends to be both a tutorial and a survey. In the first part of the paper, we present in a unified manner the stability tests expressed in terms of the Lyapunov matrix that we have developed in the past decade. For a better understanding, we organize the results according to the progression in the research: the necessary stability conditions, the infinite criterion of stability, and the finite criterion (see the note to the reader below for the paper terminology). For clarity of presentation, we consider the multiple pointwise delay case and prove the technical results in the appendix. In the second part, we give a diagnosis of the advances obtained for other classes of delay systems studied so far, highlighting the challenges raised by each class. Finally, in the last part of the paper, we discuss ongoing research on the topic and possible extensions.

More precisely, preliminaries on functionals with prescribed derivative are given in Section \ref{sec:preliminaries}, followed by fundamental stability/instability theorems for the functional in Section \ref{sec:fundamentalresults}. Section \ref{sec:necessary_stab} is devoted to the proof of necessary stability conditions and Section \ref{sec:inifinite} to the criterion of stability achieved in an infinite number of mathematical operations. In Section \ref{sec:resultsonS}, we introduce a compact set that allows the computation of constants, thus enables the presentation of a stability criterion that requires a finite number of mathematical operations  in Section \ref{sec:Finitecriterion}. We review the results obtained for other classes of delay systems in Section \ref{sec:state_of_art}. Finally, in Section \ref{sec:OngoingandFuture}, we discuss ongoing and future research, and the paper ends with conclusion.

%%%%%%%%%%%%%%%%%%%%%%%
% In the first part of the paper, we present in an organized manner the progression in the stability tests reader to the different stability criteria for linear time-delays systems. The proof strategy is explained in details in the multiple pointwise delay case. description of the organisation of the paper.
% For clarity of presentation, most of the technical results are presented in the appendices.
% In this part of the paper we have reorganized and clarified the successive steps, scattered in publications in major journals and conferences of the field of automatic control. We are the authors of the main contributions, We acknowledge the influence of the work of V. Kharitonov who set the basis for this research, the results of other researchers that were crucial to complete the results, and the work of Ph.D. who students who contributed to the develo
% The main extensions of the results to other classes of delay systems are explained, highlighting their peculiarities, and the corresponding references are given.
% In the last part of the paper, we discuss ongoing research as well as well extensions we believe are promising.

% The purpose of the present contribution is threefold. First, we describe briefly our early work and attempts to understand this interesting problem. Next, we explain in details the approach that we currently use. Then, we present a diagnosis of the advances obtained, highlighting the challenges raised by each class of systems we have studied so far. In the conclusion, we outline remaining gaps and possible future research directions.

\textit{Notation:} The space of piecewise continuous, continuous and continuously differentiable functions defined on a connected set $T\subset\R$ with values in a real finite dimensional space $X$ are represented by $\mathbf{PC}(T,X)$, $\mathbf{C}(T,X)$ and $\mathbf{C}^{(1)}(T,X)$, respectively. At the boundary points of the set $T$, one-sided continuity or the existence of a one-sided derivative are assumed. In the paper the Euclidean norm for vectors and the induced norm for matrices is denoted by $\|\cdot\|$. For functions $\ph\in\PC$ (where $H$ is a positive number), we use the uniform norm
\begin{equation*}
	\|\ph\|_\Hn=\sup_{\theta\in [-H,0]}\|\ph(\theta)\|.
\end{equation*}
The notation $A>0$ ($A \> 0$, $A \not \geqslant 0$) means that the symmetric matrix $A$ is positive definite (positive semidefinite, not positive semidefinite). The Kronecker product is denoted as $\otimes$ (it has higher priority than any other mathematical operation in this paper), the vectorization of matrix $A$ as $\mathrm{vec}(A)$, and the inverse operation as $\mathrm{vec}^{-1}$. The square block matrix with $i$-th row and $j$-th column element $A_{ij}$ is denoted by $\left[A_{ij}\right]_{i,j=1}^{r}$. The minimum eigenvalue of a matrix $A$ is represented by $\lambda_{\min}(A)$. The function that maps $y$ to the least integer greater or equal to $y$ is denoted by $\lceil y\rceil$. The identity $q\times q$ matrix is denoted by $I_q$.

\textit{Note to the reader:} In this paper, we use the following terminology: ``criterion'' (\textit{plural} criteria) means necessary and sufficient condition. \textit{Infinite criterion} means that it requires a test with an infinite number of mathematical operations, while \textit{Finite criterion} means that it can be tested by a finite number of mathematical operations.

%....I am not sure we can support this terminology with some reference\\
% Vladimir considers a criterion is NS conditions, but they are many criterion in math etc were the word is used for either necessary or sufficient, meanin "condition. That is why I thing we should put a note...\\
% e.g.: A requirement necessary for a given statement or theorem to hold. Also called a condition.\\
% SEE ALSO
% Brown's Criterion, Cauchy Criterion, Euler's Criterion, Gauss's Criterion, Korselt's Criterion, Leibniz Criterion, Pocklington's Criterion, Vandiver's Criteria, Weyl's Criterion}\\

As the system for which we present the results is linear with pointwise delays, we will use in-distinctively ``stable'' and ``exponentially stable''.

%%%%%%%%%%%%%%%%%%%%%%%%%%%%%%%%%%%%%%%%%%%%%%%%%%%%%
\section{Preliminaries} \label{sec:preliminaries}

\subsection{Basic system definitions}

Consider a linear system of the form
\begin{equation} \label{eqSyst}
	\dot{x}(t)=\sum_{j=0}^{m}A_{j}x(t-h_{j}),
\end{equation}
where $A_0$, $\ldots$, $A_m$ are real constant $n\times n$ matrices, $h_0=0$, $h_1$, $\ldots$, $h_m$ are positive delays, and $\max\{h_1,\ldots,h_m\}=H$. We additionally assume that there exists $k\in\{1,\ldots,m\}$ such that $A_k$ is nontrivial.

Without any loss of generality the initial time instant can be set equal to zero. The initial function $\ph$ is taken from the space $\PC$. The function $x(t,\ph)$, $t\>-H$, is called the solution, if it is continuous on $[0,\infty)$, satisfies~\eqref{eqSyst} almost everywhere, and $x(\theta,\ph)=\ph(\theta)$, $\theta\in[-H,0]$. The restriction of the solution to the interval $[t-H,t]$ is denoted by
\begin{equation*}
	x_t(\ph):\theta\rightarrow x(t+\theta,\ph), \;\;\theta\in[-H,0].
\end{equation*}

\begin{defn} \label{expst}
System~\eqref{eqSyst} is said to be \textit{exponentially stable}, if there exist constants $\gamma\> 1$ and $\sigma>0$, such that
\begin{equation*}
	\|x(t,\ph)\|\<\gamma e^{-\sigma t} \|\ph\|_\Hn,\;\; t\> 0.
\end{equation*}
\end{defn}

The matrix-function $K$, satisfying the equation 
\begin{equation*} \label{KprD}
	\dot{K}(t)=\sum_{j=0}^{m}A_{j}K(t-h_{j})
\end{equation*}
almost everywhere on $[0,\infty)$ with the initial conditions
\begin{equation*} \label{KprI}
	K(0)=I_n,\quad K(t)=0,\;\;t<0,
\end{equation*}
is called the \textit{fundamental matrix} of system~\eqref{eqSyst}. One can show that it satisfies also the equation 
\begin{equation*}
	\dot{K}(t)=\sum_{j=0}^{m}K(t-h_{j})A_{j},\;\; t\> 0,
\end{equation*}
almost everywhere.

The expression of the solution on $[0,\infty)$ for a given initial function $\ph$ is given by the Cauchy formula:
\begin{equation} \label{Cauchy formula}
	x(t,\ph)=K(t)\ph(0)+\sum_{j=1}^{m} \int_{-h_{j}}^{0} K(t-\theta-h_{j}) A_{j}\ph(\theta)\,d\theta.
\end{equation}

%%%%%%%%%%%%%%%%%%%%%%%%%%%%%%%%%%%%%%%%%%%%%%%%%%%%%%%%%%%%%%%%%%%%%%

\subsection{Kharitonov's approach}

We remind now the main definitions and results of the Kharitonov's approach.

In this paper we are using the right upper derivative along the solutions of system~\eqref{eqSyst}:
\begin{equation*}
	\D v(\ph) =\overline{\lim\limits_{t\to +0}} \dfrac{v(x_t(\ph))-v(\ph)}{t}.
\end{equation*}

According to \cite{KharitonovandZhabko2003}, for any positive definite matrix $W$, the functional $v_0$ satisfying
\begin{equation*}
	\D v_0(\ph) =-\ph^T(0)W\ph(0)
\end{equation*}
necessarily has the form
\begin{equation} \label{v0}
    \begin{aligned}
    	&v_{0}(\ph)=\ph^{T}(0)U(0)\ph(0) +2\ph^{T}(0) \\
    	&\cdot \sum_{j=1}^{m}\int_{-h_{j}}^{0} U(-\theta-h_j)A_{j} \ph(\theta)\,d\theta  +\sum_{i=1}^{m}\int_{-h_i}^{0} \ph^{T}(\theta_{1})A_i^T \\
    	&\cdot\sum_{j=1}^{m} \int_{-h_{j}}^{0} U(\theta_{1}+h_i-\theta_{2}-h_{j})A_{j} \ph(\theta_{2})\,d\theta_{2}\,d\theta_{1},
    \end{aligned}
\end{equation}
where the matrix-valued continuous function
\begin{equation} \label{ec:defiU}
	U(\tau)=\int_{0}^{\infty} K^{T}(t)WK(t+\tau)dt, \;\;\tau\in\R,
\end{equation}
is the \textit{delay Lyapunov matrix}. This definition requires the exponential stability assumption of the system in order to ensure convergence of the integral. This is overcome in the next definition.

\begin{defn}[\cite{Kharitonov2013}] \label{def:matriz_Lyap}
Let $W$ be a positive definite matrix. The delay Lyapunov matrix $U(\tau)$, $\tau\in\R$, of system~\eqref{eqSyst} is a function satisfying the following four properties:
\begin{enumerate}
\item[1.] Continuity property:
\begin{equation*}
    U\in\Co,
\end{equation*}
\item[2.] Dynamic property:
\begin{equation*} \label{svD}
    U^{\prime}(\tau)=\sum_{j=0}^{m}U(\tau-h_{j})A_{j}, \;\;\tau> 0,
\end{equation*}
\item[3.] Symmetry property:
\begin{equation*} \label{svS}
    U(\tau)=U^{T}(-\tau), \;\;\tau\in\R,
\end{equation*}
\item[4.] Algebraic property:
\begin{equation*} \label{svA}
    \sum_{j=0}^{m}\left(U(-h_{j})A_{j}+A_{j}^{T}U(h_{j})\right)=-W.
\end{equation*}
\end{enumerate}
\end{defn}

These properties play the role of the Lyapunov equation in the delay-free case. If the \textit{Lyapunov condition} holds (i.\,e., system~\eqref{eqSyst} has no eigenvalues that are symmetric with respect to zero), matrix $U$ exists and is unique. It can be constructed via the semi-analytic procedure when delays in system~\eqref{eqSyst} are multiple of a basic one \citep{Kharitonov2013}. In the single-delay case, the matrix $U$ is given by the following simple formula:
\begin{equation*} \label{vectorized U}
    \begin{aligned}
        &U(\tau)= \\
        &
        \begin{cases}
            \mathrm{vec}^{-1} \(
            \begin{bmatrix}
                I_{n^2} & 0
            \end{bmatrix}
            e^{L\tau }M^{-1}
            \begin{bmatrix}
                0 \\
                -\mathrm{vec}(W)
            \end{bmatrix}
            \), & \tau\>0, \\
            U^T(-\tau), & \tau\<0,
        \end{cases}
    \end{aligned}
\end{equation*}
where
\begin{equation*}
    \begin{aligned}
        L &=
        \begin{bmatrix}
            A_0^T \otimes I_{n} & A_1^T \otimes I_{n} \\
            -I_{n} \otimes A_1^T & -I_{n} \otimes A_0^T
        \end{bmatrix}, \\
        M &=
        \begin{bmatrix}
            I_{n^2} & 0 \\
            A_0^T \otimes I_{n}+I_{n} \otimes A_0^T & A_1^T \otimes I_{n}
        \end{bmatrix} \\
        &+
        \begin{bmatrix}
            0 & -I_{n^2}\\
            I_{n} \otimes A_1^T & 0
        \end{bmatrix}
        e^{LH}.
    \end{aligned}
\end{equation*}
When the construction cannot be reduced to the solution of a delay-free linear time-invariant system, one can resort to approximation methods: \cite{Jarlebringetal2011, HuescaMondie2009, Kharitonov2013, EgorovKharitonov2018, doi:10.1137/18M1209842}.

The functional of the form
\begin{equation*} \label{v integral form}
    \begin{aligned}
	    v(\ph)&=v_{0}(\ph) \\
	    &+\sum_{j=1}^{m} \int_{-h_j}^{0} \ph^{T}(\theta) \left[W_j+(\theta+h_j)W_{m+j}\right] \ph(\theta)\,d\theta,
    \end{aligned}
\end{equation*}
where $W_j$, $j=1,...,2m$, are positive definite matrices, such that
\begin{equation*}
    W_0=W-\sum_{j=1}^{m}\left(W_j+h_jW_{m+j}\right)
\end{equation*}
is also positive definite, was introduced in \cite{KharitonovandZhabko2003} for system~\eqref{eqSyst} satisfying the Lyapunov condition. Its derivative is given by
\begin{equation*}
    \begin{aligned}
		\D v(\ph)=&-\sum_{j=0}^{m} \ph^{T}(-h_j)W_j\ph(-h_j) \\ &-\sum_{j=1}^{m}\int_{-h_j}^{0}\ph^{T}(\theta)W_{m+j} \ph(\theta)\,d\theta.
	\end{aligned}
\end{equation*}
The functional $v$ allowed answering a long-standing question regarding the functional lower bound. As reported in Lemma 3.4 in \cite{Kharitonov2013}, if system~\eqref{eqSyst} is exponentially stable, there exists $\alpha>0$, such that
\begin{equation*} \label{lower bound functional complete}
	v(\ph)\>\alpha\|\ph(0)\|^2, \;\;\ph\in\PC.
\end{equation*}
This lower bound, along with the fact that the derivative of the functional along the solutions of system~\eqref{eqSyst} includes terms of the same nature as the functional itself proved to be useful in the robust stability analysis and in the determination of exponential estimates of the system response. Because of the mandatory stability assumption of the nominal system, one had to resort to well established stability results in the frequency or time domain frameworks.

These fundamental results, and the strong analogy with the delay-free linear system Lyapunov stability theory raise naturally the following question: \textit{is it possible to find stability conditions in terms of the Lyapunov matrix}?

%%%%%%%%%%%%%%%%%%%%%%%%%%%%%%%%%%%%%%%%%%%%%%%%%%%
\section{Fundamental stability theorems} \label{sec:fundamentalresults}

A first requisite for carrying out this task is presenting appropriate stability theorems. It turns out that some of the terms in functional $v$ added to $v_0$ are helpful to tackle robust stability problems but are superfluous for presenting a quadratic lower bound. Even more, they are detrimental to stability analysis. 

We introduce next the Lyapunov-Krasovskii functional $v_1$ of the form
\begin{equation*}
	v_{1}(\ph)=v_{0}(\ph) +\int_{-H}^{0}\ph^{T}(\theta) W \ph(\theta)\,d\theta,
\end{equation*}
where $v_0$ is defined in~\eqref{v0}. In the case of exponential stability,	$v_{1}$ can be rewritten as
\begin{equation} \label{v1 integral form}
	v_{1}(\ph)=\int_{-H}^{\infty}x^{T}(t,\ph)Wx(t,\ph)\,dt.
\end{equation}
 Its derivative along the solutions of system~\eqref{eqSyst} is
\begin{equation} \label{eqDv1}
	\D v_1(\ph)=-\ph^T(-H)W\ph(-H).
\end{equation}
For this functional, it is possible to formulate the following stability theorem, which is a direct consequence of Lemma 3.4 in \cite{Kharitonov2013}. Its proof is given in Appendix \ref{app:KharZabMod}.

\begin{thm} \label{KharZabMod}
If system~\eqref{eqSyst} is exponentially stable, then there exist positive numbers $\alpha_0$ and $\alpha_1$ such that for any $\ph\in\PC$
\begin{equation} \label{lower bound functional with integral}
	v_1(\ph)\>\alpha_0\|\ph(0)\|^2 +\alpha_1\int_{-H}^0\|\ph(\theta)\|^2\,d\theta.
\end{equation}
\end{thm}

Functional $v_1$ also allows proving a crucial instability result for the sufficiency, introduced in \cite{Egorov2014}. Indeed, by negation, this result is a sufficiency stability theorem. The proof is presented in Appendix \ref{app:ZhMedvMod}.

\begin{thm} \label{ZhMedvMod}
If system~\eqref{eqSyst} is unstable and satisfies the Lyapunov condition, then for every $\beta>0$ there exists a function $\widehat{\ph}\in\Cf$ such that
\begin{equation*} \label{MedvZhIn}
	v_1(\widehat\ph)\< -\beta.
\end{equation*}
\end{thm}

%%%%%%%%%%%%%%%%%%%%%%%%%%%%%%%%%%%%%%%%%%%%%%%%%%%%%%%%%%
\section{Necessary stability conditions} \label{sec:necessary_stab}

The deduction of necessary stability conditions in terms of the Lyapunov matrix from Theorem \ref{KharZabMod} relies on two key steps: the first one is the determination of the revealing initial function, and the second is the introduction of an instrumental bilinear functional.

%%%%%%%%%%%%%%%%%%%%%%%%%%%%%%%%%%%%%%%%%%%%%%%%%%%
\subsection{A suitable initial function}

It appeared in the early attempts reported in \cite{MoOchOch2011} that the substitution of simple choices of the initial function into the expression of the functional, combined with the quadratic lower bound that is satisfied when the system is stable, leads to a variety of necessary stability conditions. For example,
\begin{equation*}
	\widehat{\ph}(\theta)=
    \begin{cases}
		\mu, & \theta =0, \\
		0, & \theta\in[-H,0),
	\end{cases}
\end{equation*}
where $\mu$ is an arbitrary vector, is such that $v_{1}(\widehat{\ph}) =\mu^{T}U(0)\mu$. Then, the lower bound~\eqref{lower bound functional with integral} and the symmetry property of $U(0)$ imply that a necessary stability condition of the system is $U(0)>0$. By trial and error,  initial functions of the form $e^{A_{0}\theta}\mu$, where $\mu$ is an arbitrary non-zero vector put into evidence the dependence on the delay Lyapunov matrix in the single-delay case \citep{Mondieetal2012}. A transformation of the initial function, valid when $A_m$ is non-singular, leads to necessary conditions for the multiple delay case in the form $U(0)>0$ and
\begin{equation*} \label{sec_nec_cond}
\begin{pmatrix}
U(0) & U^{T}(\tau) & U^{T}(H) \\
U(\tau) & U(0) & U^{T}(H-\tau) \\
U(H) & U(H-\tau) & U(0)
\end{pmatrix}
\> 0,\;\; \tau\in [0,H],
\end{equation*}
reported in \cite{EgorovMondie2013b}. A decisive step was the observation that, in the single-delay case, $e^{A_{0}\theta}$ is indeed the expression of the fundamental matrix $K(\theta)$, $\theta\in [0,H]$. As a consequence, the following initial function was introduced in \cite{EgorovMondie2014} for addressing the multiple delay case:
\begin{equation} \label{initial function definition}
	\psi_r(\theta)=\sum_{i=1}^{r}K(\theta +\tau_{i})\gamma_{i}, \;\;\theta\in[-H,0],
\end{equation}
where $\tau_{i}\in[0,H]$ and $\gamma_{i}\in\R^n$, $i=\overline{1,r}$.

An important result (see, \cite{Egorovetal2017}) is that it is possible to approximate any continuously differentiable function by a function of the form~\eqref{initial function definition} with equidistant points $\tau_1$, $\tau_2$, $\ldots$, $\tau_r$.

\begin{lem} \label{FunsMix}
For any $\ph\in\Cf$ and any $\eps>0$, there exists a function $\psi_r$ of the form~\eqref{initial function definition} with
\begin{equation} \label{equidistant tau}
	\tau_i=(i-1)\delta_r,\;\;i=\overline{1,r},\quad\delta_r =\dfrac{1}{r-1}H,
\end{equation}
such that $\left\|\ph-\psi_r\right\|_\Hn<\eps$.
\end{lem}

% The proof can be found in Appendix A of \cite{Egorov2014}.

% The key role that the function $\psi_r$ plays in obtaining the necessary and sufficient stability conditions is described in more detail in the following subsections.

%%%%%%%%%%%%%%%%%%%%%%%%%%%%%%%%%%%%%%%%%%%%%%%%%%%%%%%%%
\subsection{The bilinear functional}

A crucial tool of our approach is the bilinear functional
\begin{equation*}
    z(\ph,\psi)=\dfrac{1}{4}\left(v_1(\ph+\psi) -v_1(\ph-\psi)\right),
\end{equation*}
where $\ph$, $\psi\in\PC$. Its explicit form is
\begin{equation} \label{bilinear functional}
    \begin{aligned}
	    &z(\ph,\psi)=\ph^{T}(0)U(0)\psi(0) \\
	    &+\ph^{T}(0)\sum_{j=1}^{m}\int_{-h_{j}}^{0}U(-\theta-h_j)A_{j} \psi(\theta)\,d\theta \\
        &+\sum_{j=1}^{m}\int_{-h_{j}}^{0}\ph^{T}(\theta) A_{j}^{T}U(\theta+h_{j})\,d\theta\psi(0) \\
        &+\sum_{i=1}^{m} \int_{-h_i}^{0}\ph^T(\theta_1) A_i^T \sum_{j=1}^{m}\int_{-h_{j}}^{0} U(\theta_{1}+h_i-\theta_{2}-h_{j}) \\
        &\cdot A_{j}\psi(\theta_{2})\,d\theta_{2} \,d\theta_{1} +\int_{-H}^{0}\ph^{T}(\theta)W\psi(\theta)\,d\theta.
    \end{aligned}
\end{equation}
Introducing in this bilinear functional the initial functions 
\begin{equation} \label{phipsi}
    \begin{aligned}
	    \ph(\theta)&=K(\tau_{1}+\theta)\mu, \\
	    \psi(\theta)& =K(\tau_{2}+\theta)\eta,\;\;\theta\in[-H,0],
    \end{aligned}
\end{equation}
with $\tau_{1}$, $\tau_{2}\in[0,H]$ and arbitrary vectors $\mu$, $\eta$ results in products and convolutions of the fundamental and the delay Lyapunov matrices. This interplay is governed by a number of new properties based on Definition \ref{def:matriz_Lyap}. In \cite{EgorovMondie2014}, they are proven without assuming system~\eqref{eqSyst} is stable. These properties allow a striking reduction, valid whether the system is stable or not, shown below. For the sake of clarity of presentation, the properties and the proof of this result are given in Appendix \ref{app:lemma Z reduction}.

\begin{lem} \label{lemma Z reduction}
For any $\tau_{1}$, $\tau_{2}\in[0,H]$
\begin{equation} \label{z functional}
	z\left(K(\tau_{1}+\cdot)\mu, K(\tau_{2}+\cdot)\eta\right) =\mu^{T}U(\tau_{2}-\tau_{1})\eta.
\end{equation}
\end{lem}

The following corollary is based on the equality $v_1(\ph)=z(\ph,\ph)$ and bilinearity of $z$.

\begin{cor} \label{thmV1AtPsi}
For any $\tau_1$, $\tau_2$, $\ldots$, $\tau_r\in[0,H]$
\begin{equation} \label{v1atPsi}
	v_{1}(\psi_r)=\gamma^{T} \left[\vph U(\tau_j-\tau_i)\right]_{i,j=1}^r \gamma,
\end{equation}
where $\psi_r$ is determined by~\eqref{initial function definition}, and $\gamma=\left(\gamma_1^T\,\ldots\,\gamma_r^T\right)^T$.
\end{cor}

The obtained equality reveals the meaning of the choice of function $\psi_r$. Whereas functional $v_1$ contains integrals, its value on function~\eqref{initial function definition} depends exclusively on a finite number of values of the Lyapunov matrix.

Functionals $v_1$ and $z$ satisfy upper bounds that will be useful ahead in the paper. 

\begin{lem} \label{lem:upper_bound}
There exists a number $\rho>0$ such that
\begin{equation*}
	\begin{aligned}
		|z(\ph,\psi)|&\< \rho \|\ph\|_\Hn \|\psi\|_\Hn, \\
		|v_1(\ph)|&\< \rho \|\ph\|_\Hn^2
	\end{aligned}
\end{equation*}
for any $\ph,\psi\in\PC$.
\end{lem}
A conservative, but simple estimate for $\alpha_2$ is
\begin{equation*}
    \alpha_2 \leq \nu(1+M_1)^2+HW,
\end{equation*}
where 
\begin{equation*}
    \nu=\sup_{\tau\in [0,H]}\|U(\tau)\|\hspace{3mm} \text{and}\hspace{3mm} M_1=\sum_{j=1}^{m}h_j\|A_j\|.
\end{equation*}

\subsection{Necessary stability conditions in terms of the Lyapunov matrix}

We are now ready to present necessary conditions for the exponential stability of system~\eqref{eqSyst}.

\begin{thm} \label{main result}
If system~\eqref{eqSyst} is exponentially stable, then
\begin{equation} \label{K_def}
    \left[\vph U(\tau_j-\tau_i)\right]_{i,j=1}^r>0,
\end{equation}
where $0\<\tau_1<\tau_2<\ldots<\tau_r\<H$.
\end{thm}
\begin{proof}
By Corollary \ref{thmV1AtPsi} and Theorem \ref{KharZabMod},
\begin{equation*}
    \begin{aligned}
	    v_{1}(\psi_r)&=\gamma^T\left[\vph U(\tau_j-\tau_i)\right]_{i,j=1}^r \gamma \\
	    &\>\alpha_0\|\psi_r(0)\|^2 +\alpha_1 \int_{-H}^0\|\psi_r(\theta)\|^2\,d\theta.
    \end{aligned}
\end{equation*}
It remains to show that the right-hand side is positive, if $\gamma \neq 0$. Assume by contradiction that it equals zero. As $\alpha_0>0$, $\alpha_1>0$, it implies that $\psi_r(\theta)=0$ for $\theta\in [-H,0]$ (as $\psi_r$ is right continuous), and in turn that $\gamma_i=0$ for every $i\in\{1,\hdots,r\}$.
\end{proof}

The following rough necessary stability condition is proven in \cite{EGOROV2015245}.  %appendix \ref{app:B}.

\begin{cor}\label{rough}
If system~\eqref{eqSyst} is exponentially stable, 
\begin{equation*}
    \|U(\tau)\|<\|U(0)\|, \quad \tau\in(0,\,H].
\end{equation*}
\end{cor}

\begin{rem}
Theorem \ref{main result} provides a family of necessary stability conditions, whose complexity increases with the parameter $r$ in~\eqref{K_def}. It is worth mentioning that the case $r=1$ reduces to the simplest condition $U(0)>0$, which is necessary and sufficient for the exponential stability of the delay-free system ($H=0$), and that for $r=2$ the stability criterion for the single-delay scalar equation obtained in \cite{EgorovMondie2013a} is recovered.
\end{rem}

The necessary stability conditions of Theorem \ref{main result} are illustrated by the stability maps of some examples. The continuous lines in the space of parameters represent the stability boundaries computed by using the D-subdivision method, and the isolated points correspond to the parameters values where condition~\eqref{K_def} holds.

The real numbers $\tau_i$, $i=\overline{1,r}$, of condition~\eqref{K_def} are chosen equidistant as in~\eqref{equidistant tau}.

\begin{example} \label{ex:MendezBarrios}
We consider the transfer function
\begin{equation*}
	G(s)=\dfrac{0.038}{s^4+0.1276s^3+9.3364s^2+1.1484s+3.0276}
\end{equation*}
with a controller $C(s)=k_p+k_de^{-sh_1}$ introduced in \cite{DiezMenNicMon2018}. The closed-loop system is described by a system of the form~\eqref{eqSyst} with one delay $h_1=5$ and matrices
\begin{equation*}
    \begin{aligned}
        &A_0= \\
        &
        \begin{pmatrix}
            0 & 1 & 0 & 0\\
            0 & 0 & 1 & 0\\
            0 & 0 & 0 & 1\\
            -3.0276-0.038k_p & -1.1484 & -9.3364 & -0.1276
        \end{pmatrix}, \\
        &A_1=
        \begin{pmatrix}
            0 & 0 & 0 & 0\\
            0 & 0 & 0 & 0\\
            0 & 0 & 0 & 0\\
            -0.038k_d & 0 & 0 & 0
        \end{pmatrix}.
    \end{aligned}
\end{equation*}
Figures \ref{fig:example_MB_r2} and \ref{fig:example_MB_r4} show the points in the space of parameters  $(k_p,k_d)$ where the condition of Theorem \ref{main result} holds for $r=2$ and $r=4$, respectively.   
% We depict the points where the condition of \Cref{main result} for $r=2$ and $r=4$ holds in the space of parameters $(k_p,k_d)$ in \Cref{fig:example_MB_r2,fig:example_MB_r4}, respectively. 
\begin{figure} 
	\includegraphics[height=49mm]{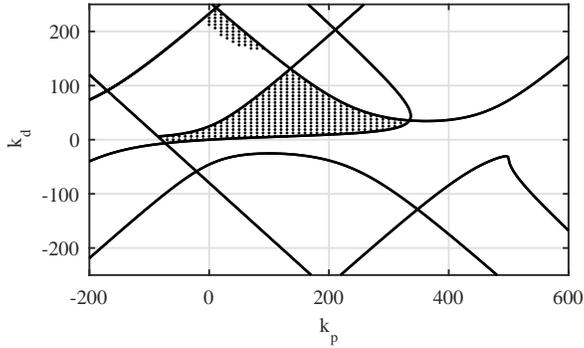}
	\caption{Example \ref{ex:MendezBarrios}, condition~\eqref{K_def} with $r=2$.}
	\label{fig:example_MB_r2}
\end{figure}
\begin{figure} 
	\includegraphics[height=49mm]{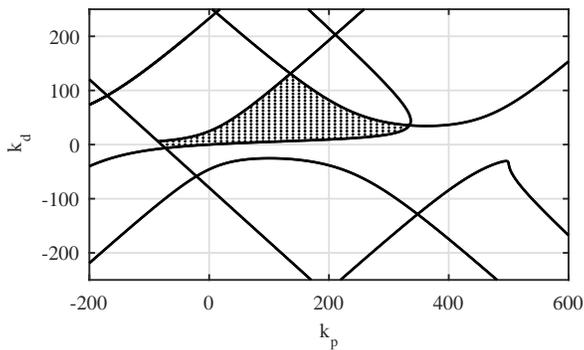}
	\caption{Example \ref{ex:MendezBarrios}, condition~\eqref{K_def} with $r=4$.}
	\label{fig:example_MB_r4}
\end{figure}
\end{example}

\begin{figure} 
	\includegraphics[height=49mm]{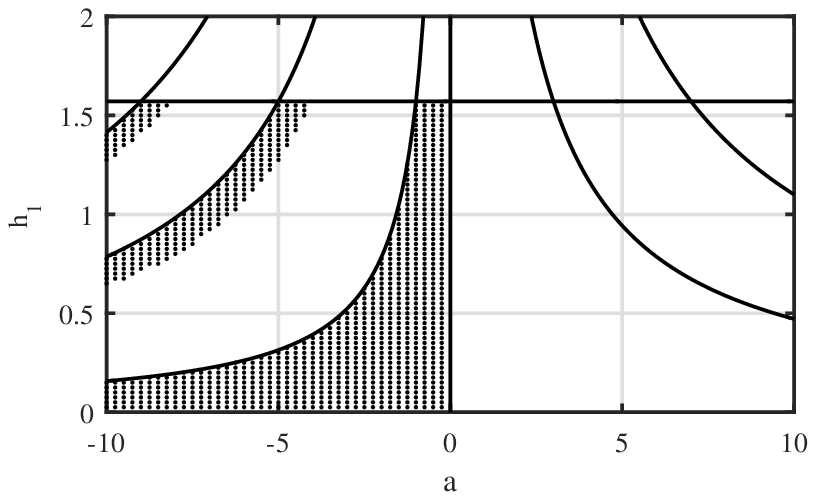}
	\caption{Example \ref{ex:matrix_A1}, condition~\eqref{K_def} with $r=3$.}
	\label{fig:example_A1_r3}
\end{figure}

\begin{figure} 
	\includegraphics[height=49mm]{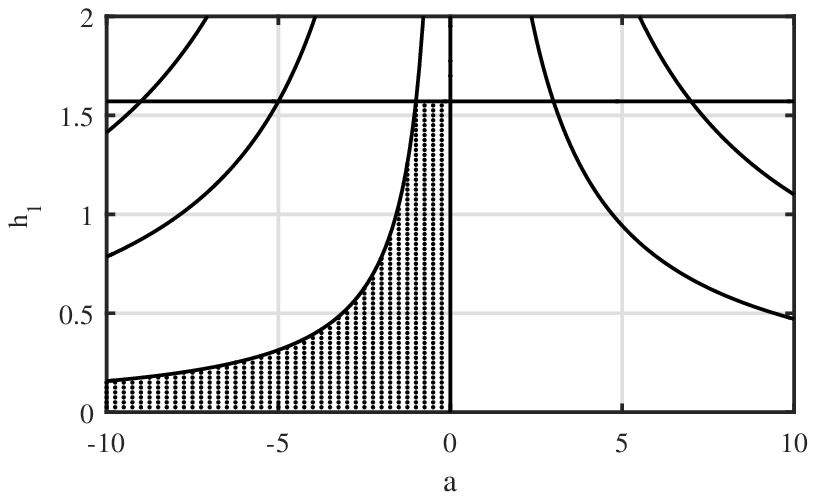}
	\caption{Example \ref{ex:matrix_A1}, condition~\eqref{K_def} with $r=6$.}
	\label{fig:example_A1_r6}
\end{figure}

\begin{example} \label{ex:matrix_A1}
We consider now the system
\begin{equation*}
\dot{x}(t)=\begin{pmatrix}
-1 & 0.5\\
0 & a
\end{pmatrix}x(t-h_1).
\end{equation*}
The points where condition~\eqref{K_def} for $r=3$ and $r=6$ holds in the space of parameters $(a,h_1)$ are depicted in Figures \ref{fig:example_A1_r3} and  \ref{fig:example_A1_r6}.
\end{example}

\begin{rem}
The examples indicate that the conservatism of the necessary stability condition~\eqref{K_def} can be reduced by increasing the number $r$. Indeed, the exact stability regions in Examples \ref{ex:MendezBarrios} and \ref{ex:matrix_A1} are achieved for $r=4$ and $r=6$, respectively. Notice that for $r=2^\ell$, $\ell=1,2,3,...$, the conditions form a hierarchy.
\end{rem}

The obtained necessary stability conditions in terms of the delay Lyapunov matrix, and their efficacy in providing accurate estimates of the stability region observed in the above examples (see more examples in \cite{Cuvasetal2017}) lead to the following question: is it possible to prove sufficiency of the conditions of Theorem \ref{main result}?

\section{Infinite criterion of stability} \label{sec:inifinite}

A first answer to the question on sufficiency is an infinite criterion of stability for systems with pointwise delays. 
% and distributed delays, is proven in full detail in \cite{Egorov2014,Egorovetal2017}.

Let us introduce matrices
\begin{equation*}
    \begin{aligned}
        \K_1&=U(0), \\
        \K_2&=
        \begin{pmatrix}
            U(0) & U(H) \\
            \ast & U(0)
        \end{pmatrix}, \\
        \K_3&=
        \begin{pmatrix}
            U(0) & U(H/2) & U(H) \\
            \ast & U(0) & U(H/2) \\
            \ast & \ast & U(0)
        \end{pmatrix}, \\
        \K_4&=
        \begin{pmatrix}
            U(0) & U(H/3) & U(2H/3) & U(H) \\
            \ast & U(0) & U(H/3) & U(2H/3) \\
            \ast & \ast & U(0) & U(H/3)\\
            \ast & \ast & \ast & U(0)
        \end{pmatrix},
    \end{aligned}
\end{equation*}
and so on. In general, for $r=2,3,\ldots$
\begin{equation*}
    \K_r =\left[\vph U\left(\dfrac{j-i}{r-1}H\right)\right]_{i,j=1}^r.
\end{equation*}

The next theorem states that the necessary stability condition of Theorem \ref{main result} is also sufficient for large enough $r$.

\begin{thm} \label{MainTheorem}
System~\eqref{eqSyst} is exponentially stable if and only if the Lyapunov condition holds and for every natural number $r$
\begin{equation} \label{UCNF}
	\K_r>0.
\end{equation}
Moreover, if the Lyapunov condition holds and system~\eqref{eqSyst} is unstable, then there exists a natural number $r$ such that
\begin{equation*}
	\K_r\not\geqslant 0.
\end{equation*}
\end{thm}
\begin{proof}
\textit{Necessity:} The system is assumed to be exponentially stable, hence, the Lyapunov condition holds. For every natural $r$, in the special case of equidistant points~\eqref{equidistant tau} condition~\eqref{K_def} in Theorem \ref{main result} takes the form~\eqref{UCNF}.

\textit{Sufficiency:} Assume that the Lyapunov condition holds, then the delay Lyapunov matrix $U$ exists. To prove the sufficiency it is enough to show that for unstable systems there exists a natural number $r$, such that
\begin{equation} \label{ContInProof0}
	\K_r\not\geqslant 0.
\end{equation}
In view of equation~\eqref{v1atPsi}, which was established without stability assumption, inequality~\eqref{ContInProof0} is equivalent to the existence of a vector $\gamma$, such that
\begin{equation} \label{ContInProof}
	\gamma^{T}\K_r\gamma <0.
\end{equation}
Suppose that system~\eqref{eqSyst} is unstable. Set an arbitrary number $\beta>0$. By Theorem \ref{ZhMedvMod}, there exists $\widehat{\ph}\in\Cf$, such that $v_1(\widehat{\ph})\<-\beta$. Functional $v_1$ is continuous at each point, i.\,e., there exists a number $\Delta>0$, such that
\begin{equation*}
	\|\widehat{\ph}-\ph\|_\Hn<\Delta \quad \Rightarrow\quad |v_{1}(\widehat{\ph})-v_{1}(\ph)| <\beta.
\end{equation*}
By Lemma \ref{FunsMix}, there exists $\psi_r$ of the form~\eqref{initial function definition} with equidistant $\tau_i$, $i=\overline{1,r}$, like in~\eqref{equidistant tau}, such that
\begin{equation*}
	\|\widehat{\ph}-\psi_r\|_\Hn<\Delta.
\end{equation*}
Hence,
\begin{equation*}
	v_{1}(\psi_r)<v_{1}(\widehat{\ph})+\beta\< 0.
\end{equation*}
As function $\psi_r$ has the form~\eqref{initial function definition}, equality~\eqref{v1atPsi} implies that
\begin{equation*}
	v_{1}(\psi_r)=\gamma^{T}\K_r\gamma <0.
\end{equation*}
Thus, inequality~\eqref{ContInProof}, equivalently~\eqref{ContInProof0}, is shown, and sufficiency is proven.
\end{proof}
% \textcolor{red}{The results and proof this section should be written as much as possible in a similar manner as those of section 7}

In order to use Theorem \ref{MainTheorem} for determining the stability of the system, one requires an \emph{infinite} number of mathematical operations, as condition~\eqref{UCNF} demands to be tested for every natural $r$.  Although the theorem states that for every unstable system there exists a number $r$ such that condition~\eqref{UCNF} does not hold, it does not provide any estimate of such number.

%%%%%%%%%%%%%%%%%%%%%%%%%%%%%%%%%%%%%%%%%%%%%%%%%%
\section{Refined results on a compact set} \label{sec:resultsonS}

The difficulty for estimating the number $r$ for which condition~\eqref{UCNF} does not hold is that the approximation error depends on $\ph$. To remedy this issue, inspired by \cite{MedveZhabko2015}, the stability theorems~\ref{KharZabMod} and~\ref{ZhMedvMod}, and the approximation error are revisited in the framework of the following relatively compact set:
\begin{equation*}
	\begin{aligned}
		\S=\left\{\vphS\right.\ph\in\Cf:& \\
        \|\ph\|_\Hn=\|\ph(0)\|&=1, \|\ph'\|_\Hn\< M\left.\vphS\right\},
	\end{aligned}
\end{equation*}
where $M=\sum_{j=0}^m\|A_j\|$.

%%%%%%%%%%%%%%%%%%% approximation %%%%%%%%%%%%%%%%%%%%%%%%%
The next lemma shows that every function from the set $\S$ can be approximated by a function $\psi_r$ constructed as in Lemma \ref{FunsMix}. In contrast with Lemma \ref{FunsMix}, an approximation error bound is now provided. The proof can be found in \cite{Egorov2016}.

\begin{lem} \label{lem:approximation_if}
For every $\ph\in\S$, there exists a function $\psi_r$ of the form~\eqref{initial function definition} with equidistant $\tau_i$, $i=\overline{1,r}$, like in~\eqref{equidistant tau}, such that $\ph(0)=\psi_r(0)$ and $\|\ph-\psi_r\|_\Hn\<\eps_r$, where
\begin{equation*}
	\eps_r=\dfrac{(M+L)e^{LH}}{1/\delta_r+L}.
\end{equation*}
Here, $L$ is such that $\|K'(t)\|\< L$, $t\in[0,H]$, almost everywhere.
\end{lem}

%\begin{lem} \label{lem:spectral_abscissa}
%	The unstable characteristic roots of system~\eqref{ec:sys} fulfills
%	\begin{equation}
%	\mathbf{Re}(s)<\hat\alpha.
%	\end{equation}
%	Here, $\hat \alpha$ is a positive number that satisfies
%	\begin{equation} \label{ec:condition_alpha_hat}
%	e^{2\hat \alpha h_1} >\dfrac{\sum\limits_{i=1}^m \|PA_i\|^2}{m\lambda_{\min}^2(Q)},
%	\end{equation}
%	where $P$ is a positive definite matrix, such that
%	\begin{equation} \label{ec:Lyapunov_equation_alpha}
%	(A_0-\hat \alpha I)^TP+P(A_0-\hat \alpha I)=-2mQ.
%	\end{equation}
%\end{lem}

%Number $\hat \alpha$ can be computed as follows \citep{TissirandHmamed1996}:
%\begin{enumerate}[Step 1:]
%	\item Set $\hat \alpha=0$ and $\Delta \hat \alpha$.
	
%	\item Solve the Lyapunov equation~\eqref{ec:Lyapunov_equation_alpha}.
	
%	\item Compute
%	\begin{equation}
%	\hat\alpha_{\min}=\dfrac{1}{2h_1} \,\ln\dfrac{\sum\limits_{i=1}^m \|PA_i\|^2}{m\lambda_{\min}^2(Q)}.
%	\end{equation}
%	If $\hat\alpha_{\min}\< \hat\alpha$, stop the process, otherwise set $\hat\alpha=\hat\alpha+\Delta \hat \alpha$ and go to Step 2.
%\end{enumerate}

The following result proven in Appendix \ref{app:theo:stability_v1_set}, strengthens Theorem \ref{KharZabMod}.

\begin{thm} \label{theo:stability_v1_set}
If system~\eqref{eqSyst} is exponentially stable, then for any $\ph\in\S$
\begin{equation} \label{lower bound functional}
	v_{1}(\ph)\>\alpha_0^{\star},
\end{equation}
where
\begin{gather*}
    \alpha_0^{\star}=-\dfrac{1}{(m+1)\lambda_{\min}(P)}, \\
    P=I_{m+1}\otimes W^{-1} \cdot \left(e_1\otimes I_n\cdot A + A^T\cdot e_1^T\otimes I_n\right), \\
    e_1=(1,0,\ldots,0)^T\in\R^{m+1},
\end{gather*}
and block-matrix
\begin{equation*}
A=\left(A_0,A_1,\ldots,A_m\right)\in\R^{n\times n(m+1)}.
\end{equation*}
\end{thm}

\begin{rem}
It is worth noting that all eigenvalues of matrix $P$ are real and $\lambda_{\min}(P)<0$ for any system whether it is stable or not. Thus, number $\alpha_0^{\star}$ can be computed for any system.
\end{rem}

The following theorem is similar to Theorem \ref{ZhMedvMod} except for the fact that it provides a \textit{computable} bound of the functional $v_1$. This bound is, indeed, the cornerstone of the stability criterion presented in the next section. The proof is given in Appendix \ref{app:instability_theorem}.

\begin{thm} \label{theo:instability}
If system~\eqref{eqSyst} is unstable and satisfies the Lyapunov condition, there exists $\widehat\ph\in\S$ such that
\begin{equation*}
	v_1(\widehat\ph)\< -\beta^{\star}
\end{equation*}
with
\begin{equation*}
	\beta^{\star} =\dfrac{\lambda_{\min}(W)}{4a} e^{-2a H} \cos^2(b),
\end{equation*}
where $a>0$ is an upper estimate of the real part of the rightmost root of system~\eqref{eqSyst}\footnote{A rough estimate is given by $M=\sum_{j=0}^m\|A_j\|$. In \cite{Gomezetal2019} an algorithm to improve this estimate is provided.} and $b$ is a unique solution of the equation
\begin{equation} \label{ec:b}
	((a H)^2+b^2)\sin^4(b)=(a H)^2
\end{equation}
on $\left(0,\frac{\pi}{2}\right)$.
\end{thm}

%Revisiting \Cref{KharZabMod,ZhMedvMod} when the initial functions belong to the set $\S$ introduced in the previous section gives a stability criterion for functional $v_1$:

%%%%%%%%%%%%%%%%%%%%%%%%%%%%%%%%%%%%%%%%%%%%%%%%%%%%%%

\section{Finite  criteria of stability} \label{sec:Finitecriterion}

The results of the previous section allow presenting stability criteria that require a finite number of mathematical operations. The bound on $r$ for which the instability condition holds is found, implying sufficiency of the conditions. 

\subsection{Finite criterion in terms of the Lyapunov matrix}

The following criterion, presented in \cite{Gomezetal2019}, depends exclusively on the Lyapunov matrix.

\begin{thm} \label{theo:NSC_condition_LM}
System~\eqref{eqSyst} is exponentially stable if and only if the Lyapunov condition holds and
\begin{equation} \label{ec:stability_criteria_LM}
	\K_r>0,
\end{equation}
where
\begin{equation} \label{ec:number_r2}
	r=1+\left\lceil He^{LH}\left(M+L\right) \left(\alpha +\sqrt{\alpha(\alpha+1)}\right) -HL\right\rceil
\end{equation}
with $\alpha=\dfrac{\rho}{\beta^{\star}}$. Here $\rho$ and $\beta^{\star}$ are determined by Lemma \ref{lem:upper_bound} and Theorem \ref{theo:instability}, respectively.
\end{thm}
\begin{proof}
The necessity was proven in Theorem \ref{MainTheorem}. In order to prove the sufficiency, we assume by contradiction that system~\eqref{eqSyst} is unstable but that the condition in~\eqref{ec:stability_criteria_LM} and the Lyapunov condition hold, thus, the existence and uniqueness of the delay Lyapunov matrix is guaranteed.

By Theorem \ref{theo:instability}, there exists $\widehat\ph\in\S$ such that $v_1(\widehat\ph) \<-\beta^{\star}$. Construct a corresponding $\psi_r$ from Lemma \ref{lem:approximation_if}.

By defining $E_r=\widehat\ph-\psi_r$, we get  
\begin{equation*}
	 v_1(\psi_r) = v_1(\widehat\ph-E_r)=z(\widehat\ph-E_r, \widehat\ph-E_r).
\end{equation*}
The bilinearity of $z$ and Lemma \ref{lem:upper_bound} imply that
\begin{equation*}
    \begin{aligned}
	    v_1(\psi_r)&= v_1(\widehat\ph)-2z(\widehat\ph,E_r)+v_1(E_r) \\
	    &\< -\beta^{\star} +2\rho\|E_r\|_\Hn+\rho\|E_r\|_\Hn^2 \\
	    &=\beta^{\star}\left(-1 +2\alpha\|E_r\|_\Hn+\alpha\|E_r\|_\Hn^2\right).
	\end{aligned}
\end{equation*}
By using Lemma \ref{lem:approximation_if} and considering the number $r$ given by~\eqref{ec:number_r2}, we have
\begin{equation*}
    \|E_r\|_\Hn\< \dfrac{(M+L)e^{L H}}{1/\delta_r+L} \< \dfrac{1}{\alpha +\sqrt{\alpha(\alpha+1)}},
\end{equation*}
which implies that
\begin{equation*}
	-1 +2\alpha\|E_r\|_\Hn +\alpha\|E_r\|_\Hn^2\< 0.
\end{equation*}
Finally, from the previous inequality and formula~\eqref{v1atPsi}, we get
\begin{equation*}
	v_1(\psi_r)=\gamma^T\K_r\gamma\< 0,
\end{equation*}
which contradicts the assumption that~\eqref{ec:stability_criteria_LM} holds.
\end{proof}

Let us revisit Example \ref{ex:matrix_A1} to illustrate the results of Theorem \ref{theo:NSC_condition_LM}.

\begin{example}[Example \ref{ex:matrix_A1} revisited] \label{ex:ex2revi}
Table \ref{tab:comp_r} displays values of $r$ from~\eqref{ec:number_r2} in the second column for some parameters $(a,h_1)$ of the system in Example \ref{ex:matrix_A1}\footnote{The numerical computation were performed in a computer with a processor  Intel Core i5 with two cores of 2.5GHz and RAM memory of 8GB.}.
\begin{table*} 
	\caption{Computation of $r$ from~\eqref{ec:number_r2} and \eqref{ec:number_r} for some parameters $(a,h_1)$ of the system in Example \ref{ex:matrix_A1}. The computation time in the third and sixth columns includes the computation of $r$, the construction of $U$, $\K_r$ and $\P_r$, and the positive definiteness test.}
	\begin{center}
		\begin{tabular}{ |c | c|c|c||c|c|c| }\hline
			Parameters $(a,h_1)$ & $ r$ (Th. 7) & [sec] (Th. 7) & Test result (Th 7) & $ r$ (Th. 8) & [sec] (Th. 8) & Test result (Th 8) \\ \hline
			$(-1.25,0.5)$ & $89$ &$1.30$ & $\K_{89}>0$ & $14$ & $1.46$ & $\K_{14}-\alpha_0\P_{14}^T\P_{14}>0$ \\
			$(-1.25,0.75)$ & $395$ &$1.57$ & $\K_{395}>0$ & $45$ & $1.36$ & $\K_{45}-\alpha_0\P_{45}^T\P_{45}>0$ \\ 
			$(1.25,0.5)$ & $79$ & $1.35$& $\K_{79}\not\geqslant 0$ & $13$ & $1.36$ & $\K_{13}-\alpha_0\P_{13}^T\P_{13}\not\geqslant0$ \\ 
			$(1.25,1.25)$ & $3416$ & $23.84$ &$\K_{3416}\not\geqslant 0$ & $257$ & $1.38$ & $\K_{257}-\alpha_0\P_{257}^T\P_{257}\not\geqslant 0$ \\
			\hline
		\end{tabular}
	\end{center}
	\label{tab:comp_r}
\end{table*}

The result of the stability test is shown in the fourth column. According to the stability criterion of Theorem \ref{theo:NSC_condition_LM}, the parameters in the first and second rows correspond to stable systems, whereas the third and fourth rows to unstable systems. The computation time in the third column includes the computation of $r$ given by ~\eqref{ec:number_r2}, the construction of $U$ and $\K_r$ and the positive definiteness test of the latter matrix. It is important to mention that constructing $\K_r$ is the most computation time consuming task.
\end{example}

The above example shows that the number $r$ of Theorem \ref{theo:NSC_condition_LM} widely differs from the number $r=6$ for which the exact stability zone is found in Example \ref{ex:matrix_A1}. The same is observed in several illustrative examples where the stability region is reached with a small number $r$; see, e.\:g., \cite{Cuvasetal2017} and the references therein. While it is true that the number $r$ from~\eqref{ec:number_r2} can be very large, a fact of theoretical significance is that by using Theorem \ref{theo:NSC_condition_LM} one can determine the stability of the system with a finite number of mathematical operations. 

\subsection{Improved stability criterion}
A current challenge is to reduce the estimate $r$ presented in Theorem \ref{theo:NSC_condition_LM}. A finite criteria of stability, established with the help of Theorem \ref{theo:stability_v1_set}, presented in \cite{Egorov2016}, reduced this estimate. A further reduction was achieved in \cite{Egorov2020}, which slashed the number $r$ by approximately two. 
 
% \textcolor{red}{AAAAAAAAAAAAAAAAAAAAAAAAAAAAAA\\
% A discussion? in the line of CDC 2020 (qualitative)\\
% MAYBE WE SHOULD ONLY DISCUSS THESE THEOREMS and give the general above description...the interested reader can go study the papers?\\
% AAAAAAAAAAAAAAAAAAAAAAAAAAAAAAAA}

In this section, we give one possible improvement, which has not been published before in the present form.

Consider the block-matrix
\begin{equation*}
	\P_r=\left(I_n,K(\delta_r),\ldots,K((r-1)\delta_r)\right) \in\R^{n\times n r}.
\end{equation*}

\begin{thm} \label{theo:finite_criterion}
System~\eqref{eqSyst} is exponentially stable if and only if the Lyapunov condition holds and
\begin{equation*}
	\K_r-\alpha_0\P_r^T\P_r>0,	
\end{equation*}
where
\begin{equation} \label{ec:number_r}
	r=1+\left\lceil He^{LH}(M+L) \left(\alpha+\sqrt{\alpha(\alpha+1)}\right)-HL \right\rceil
\end{equation}
with $\alpha=\dfrac{\rho}{\beta^{\star}+\alpha_0}$ and $\alpha_0\in(0,\alpha_0^{\star})$. Here $\rho$, $\alpha_0^{\star}$ and $\beta^{\star}$ are determined by Lemma \ref{lem:upper_bound}, Theorem \ref{theo:stability_v1_set} and \ref{theo:instability}, respectively.
\end{thm}
\begin{proof}
\textit{Necessity:} By using the lower bound~\eqref{lower bound functional} and equation~\eqref{v1atPsi}, we have, for every $\gamma\in\R^{nr}$ such that $\psi_r(0)\neq 0$,
\begin{equation*}
	\begin{aligned}
		\gamma^T\left(\K_r -\alpha_0\P_r^T\P_r\right)\gamma
        &=v_1(\psi_{r})-\alpha_0\|\psi_{r}(0)\|^2 \\
        &> v_1(\psi_{r})-\alpha_0^{\star}\|\psi_{r}(0)\|^2 \>0.
	\end{aligned}
\end{equation*}
For the case $\psi_r(0)=0$, $\gamma\neq0$, the inequality
$$
\gamma^T\left(\K_r -\alpha_0\P_r^T\P_r\right)\gamma>0
$$
remains valid by Theorem \ref{MainTheorem}.

\textit{Sufficiency:} As in the proof of Theorem \ref{theo:NSC_condition_LM} we obtain
\begin{equation*} \label{ec:SC_eqfunct}
    \begin{aligned}
	    v_1(\psi_r)&=v_1(\widehat\ph)-2z(\widehat\ph,E_r)+v_1(E_r) \\
	    &\< -\beta^{\star} +2\rho \|E_r\|_\Hn+\rho\|E_r\|_\Hn^2.
    \end{aligned}
\end{equation*}
Since $\|\psi_r(0)\|=\|\ph(0)\|=1$,
\begin{equation*}
    \begin{aligned}
	    \gamma^T&\left(\K_r-\alpha_0\P_r^T\P_r\right)\gamma =v_1(\psi_r) -\alpha_0\|\psi_r(0)\|^2 \\
	    &\< -\left(\alpha_0+\beta^{\star}\right) +2\rho \|E_r\|_\Hn+\rho\|E_r\|_\Hn^2 \\
	    &\< \left(\alpha_0+\beta^{\star}\right)\left(-1 +2\alpha \|E_r\|_\Hn+\alpha\|E_r\|_\Hn^2\right).
    \end{aligned}
\end{equation*}
The rest of the proof is similar to the one of Theorem \ref{theo:NSC_condition_LM}.
\end{proof}

Computed numbers $r$ in the stability criterion of Theorem \ref{theo:finite_criterion}, stability test result and computation time for the system parameters $(a,h_1)$ of Example \ref{ex:ex2revi} are shown Table \ref{tab:comp_r}. The reduction of the estimate order is evident. Notice, however, that in Theorem \ref{theo:finite_criterion} the stability criterion does not depend uniquely on the Lyapunov matrix, but also on the fundamental matrix of the system. 

\section{State of the art} \label{sec:state_of_art}

We have been able to present stability tests in terms of the Lyapunov matrix of several clases of delay systems. We outline in this section the progress made up to now.

% An interesting fact is that, for ordinary differential systems, the stability condition preserves the form obtained for the retarded type case. Of course, what is different is the underlying delay Lyapunov matrix, which is constructed with the help of the three properties~\eqref{svD},~\eqref{svS} and~\eqref{svA} that are defined by the class of system.

%%%%%%%%%%%%%%%%%%%%%%%%%%%%%%%%%%%%%%%%%%%%%%%%%%%%%%%%%
\subsection{Systems with pointwise and distributed delays}

The approach described in this paper was extended  to the case of systems with multiple pointwise and distributed delays of the form
\begin{equation*} \label{sys}
	\dot{x}(t)=\sum_{j=0}^{m}A_{j}x(t-h_{j}) +\int_{-H}^{0}G(\theta) x(t+\theta)\,d\theta,
\end{equation*}
where $G$ is a real piecewise continuous matrix function defined on $[-H,0]$.

The necessary conditions were presented in \cite{CuvasMondie2016}  and the infinite criterion of stability is proven in \cite{Egorovetal2017}. The finite criterion is available in \cite{Castano2022}. The form~\eqref{K_def} of the necessary conditions is the same; of course the Lyapunov matrix is now the one corresponding to the distributed delay system. For these systems, a semi-analytic construction of the delay Lyapunov matrix is possible when, on the one end, the delays and the discontinuity points of $G$ are commensurate, and on the other hand, the kernel can be expressed as the solution of a linear delay-free system on every segment of discontinuity. The basic construction method introduced in \cite{KHARITONOV2006610} is clarified in \cite{Aliseyko2019}. Application examples with constant, polynomial, exponential and gamma distributions are reported in \cite{CUVAS2015239, Aliseyko2019, 8619265, ALISEYKO201919,JuarezMonKhar2020}.

%%%%%%%%%%%%%%%%%%%%%%%%%%%%%%%%%%%%%%%%%%%%%%%%%%%%%%%%%
\subsection{Systems of neutral type}

We have obtained necessary stability conditions, also of the form~\eqref{K_def}, in \cite{Gomezetal2016b} for systems of neutral type described by
\begin{equation*} \label{ec:sys_neutral}
	\dfrac{d}{dt}\left(\vph x(t)-Dx(t-H)\right)=A_0x(t)+A_1x(t-H),
\end{equation*}
where $D$, $A_0$ and $A_1$ belong to $\R^{n\times n}$, and $H>0$ is the delay. The neutral case is, as usual, more tricky, and special care has to be paid to the jump discontinuities of the system. The generalization of the necessary conditions to the multiple commensurate delays case is carried out in \cite{Gomezetal2017a}. As for retarded type systems, the attainment of the stability criterion for neutral type systems is based on the introduction of the compact set
\begin{equation*}
	\begin{aligned}
		\S_{\mu}=\left\{\vphS\right.\ph\in\Cf:& \\
        \|\ph\|_\Hn=\|\ph(0)\|&=1, \|\ph'\|_\Hn\< \mu M\left.\vphS\right\},
	\end{aligned}
\end{equation*}
where
\begin{equation*}
    M=\|A_0\|+\|A_1\|\:\:\text{and}\;\;\mu=\dfrac{d}{1-\rho}
\end{equation*}
with $d\>1$ and $\rho\in(0,1)$ such that $\|D^k\|\< d\rho^k$, $k=0,1,2,\ldots$. Some advances on the finite stability criterion of neutral type delay system were presented in \cite{Gomezetal2018} when matrix $D$ is such that $\|D\|<1$. The general result in the form of Theorem \ref{theo:NSC_condition_LM} is proven in \cite{GomezEgMon2021}.
Again, the effective testing of the conditions depends on the availability of the delay Lyapunov matrix of the system. For the construction method introduced in \cite{Kharitonov2013} and examples, the reader is refereed to \cite{2016GomCuvCDC2016}.

% \begin{thm}\cite{Gomezetal2018}
% Assume that the matrix $D$ satisfies $\|D\|<1$. System~\eqref{ec:sys_neutral} is exponentially stable if and only if the Lyapunov condition holds and
% \begin{equation}
% 	\K_r-\alpha_1\A_r>0,
% \end{equation}
% where
% \begin{equation}
% 	r=...,
% \end{equation}
% \begin{equation}
% 	\A_r=...
% \end{equation}
% \end{thm}

%%%%%%%%%%%%%%%%%%%%%%%%%%%%%%%%%%%%%%%%%%%%%%%%%%%%%%%%
\subsection{Linear periodic systems with delays}

In \cite{GomezOchoaMondie2016}, inspired by the developments of converse results presented in \cite{LetyaginaZhabko2009}, we provide necessary stability conditions for the class of systems described by
\begin{equation} \label{ec:sys_periodic}
	\dot{x}(t)=\sum_{j=0}^{m}A_{j}(t)x(t-h_{j}),
\end{equation}
where $A_j$, $j=\overline{0,m}$, are matrices of continuous coefficients with period $T$, i.\,e., $A_j(t)=A_j(t+T)$, $t\>0$, $j=\overline{0,m}$, with range in $\R^{n\times n}$, and the delays $h_0$, $h_1$, $\ldots$, $h_m$ are constant. The stability conditions are similar to those obtained in the time-invariant case.

\begin{thm}
If system~\eqref{ec:sys_periodic} is exponentially stable, then the following condition holds for any positive integer $r$, any $\tau_i\in[0,H]$, $i=\overline{1,r}$, and any $t\in[0,T)$:
\begin{equation*}
    \left[U(t-\tau_i,t-\tau_j)\vph\right]_{i,j=1}^r\> 0.
\end{equation*}
\end{thm}

The delay Lyapunov matrix for system~\eqref{ec:sys_periodic} is a function of two arguments. In the single-delay case with $T=h_1=H$ it is solution of the partial differential equation system
\begin{equation*} \label{sysPartial}
	\begin{aligned}
	    \dfrac{\partial U(\tau_{1},\tau_{2})}{\partial\tau_{1}} &=-A_{0}^T(\tau_{1})U(\tau_{1},\tau_{2})-A_1^T(\tau_1)V(\tau_{1},\tau_{2}),\\
	    \dfrac{\partial V(\tau_{1},\tau_{2})}{\partial\tau_{2}} &=-V(\tau_{1},\tau_{2})A_{0}(\tau_{2}) -U(\tau_{1},\tau_{2})A_{1}(\tau_{2}),
	    \end{aligned}
\end{equation*}
subject to the boundary conditions
\begin{equation*} \label{initialCond}
	\begin{aligned}
    	\dfrac{d}{d\tau}\left(\vphS U(\tau,\tau)\right) &=-A_0^T(\tau)U(\tau,\tau) -A_1^T(\tau)V(\tau,\tau) \\
	    &-U(\tau,\tau)A_0(\tau) -U(\tau-T,\tau)A_1(\tau)-W,\\
        V(\tau-T,\tau)&=U(\tau,\tau),\\
	    U(\tau-T,0)&=U(\tau,T),\\
	    V(\tau-T,0)&=V(\tau,T),
    \end{aligned}
\end{equation*}
where $\tau_1$ and $\tau_2$ are such that $\tau_2-\tau_1\in[0,T]$, $\tau_2\in[0,T]$, and $\tau\in[0,T]$. 

The computation of the delay Lyapunov matrix is a challenging issue that is addressed in \cite{GOMEZEgMonZhab2019,MichielGomez2020}.

%%%%%%%%%%%%%%%%%%%%%%%%%%%%%%%%%%%%%%%%%%%%%%%%%%%%%%%%
\subsection{Difference equations in continuous time and integral delay equations}

We have extended the present approach to some classes of non-differential linear systems, namely, difference equations of the form
\begin{equation*} \label{eq:sys}
	x(t)=\sum_{j=1}^{m}{A_{j}x(t-h_{j})},
\end{equation*}
where $x(t)\in\R^n$, $A_1$, $\ldots$, $A_m$ are constant real $n\times n$ matrices, and $h_1$, $h_2$, $\ldots$, $h_m$ are the delays, $\max\{h_1,\ldots,h_m\}=H$, and also to integral equations of the form
\begin{equation*}
	x(t)=\int_{-H}^{0}G(\theta)x(t+\theta)\,d\theta
\end{equation*}
with piecewise continuous kernel $G$, and delay $H$. These equations appear in the difference operator of neutral type systems, and in the modeling of distributed parameter systems. The fundamental matrix and Cauchy formula of integral equations is reported in \cite{MelchorKhar2010}. The definition of the Lyapunov matrix,
\begin{equation*} \label{U_tau}
    U(\tau)=\int_{0}^{\infty}\left(K(t)-K_0\right)^T WK(t+\tau)dt,
\end{equation*}
insures it is well-defined for all $\tau\in\R$ when the system is stable. Here,
\begin{equation*} \label{MtxFunCI1}
    K_0 =\left(\sum_{j=1}^{m}{A_{j}} -I_n\right)^{-1}
\end{equation*}
for difference equations, and
\begin{equation*} \label{MtxFunCI}
    K_0 =\left(\int_{-h}^0 G(\theta) d\theta -I_n\right)^{-1}
\end{equation*}
for integral equations. In both cases, the necessary conditions have the same general form
\begin{equation} \label{eq:kkh}
	\left[F(\tau_i,\tau_j)\right]_{i,j=1}^{r}>0,
\end{equation}
where $0\<\tau_1<\tau_2<\ldots<\tau_r\<H$,
\begin{equation*} \label{eq:fff}
	F(\tau_i,\tau_j)=U(0)-U(\tau_j)-U(-\tau_i)+U(\tau_j-\tau_i).
\end{equation*}
Necessary stability conditions for difference equations in continuous time are presented in \cite{Rochaetal2016}. Necessary conditions for integral delay equations found in \cite{delValleetal2018} are shown to be sufficient in \cite{Ortiz2022}.

% The Dini upper right-hand derivative of a functional was introduced to address the jump discontinuities inherent to this class of systems.

Of course, the Lyapunov matrix in~\eqref{eq:kkh} depends on the class of system under consideration. For integral equations, the construction of $U$ for the case of kernels that are solution of a system of ordinary differential equations is available in \cite{ORTIZ201991}. For difference equations in continuous time, the construction method is presented in \cite{ROCHA20176507}: for the commensurate delay case, it is trivial, while for the incommensurate case, one has to resort to approximations.

% \begin{rem}
% It may be possible to prove conditions in the simpler form~\eqref{K_def} obtained for differential systems by using another definition of the matrix $U$ and different initial functions in the proof of the result.
% \end{rem}

%%%%%%%%%%%%%%%%%%%%%%%%%%%%%%%%%%%%%%%%%%%%%%%%%%%%%%%%%
\subsection{Linear systems with incommensurate delays}

The general machinery employed in the paper remains valid in the incommensurate delay case. However, when the delays in~\eqref{eqSyst} are not commensurate, the semi-analytic construction for matrix $U$ is not possible. Thus, only an approximation is accessible. In \cite{EgorovKharitonov2018}, a computable bound on the approximation error of the Lyapunov matrix of a stable incommensurate delay system and a nearby stable commensurate one is presented, and it is shown how the approximate functional allows to analyze stability. These results allowed the presentation of necessary stability conditions for the incommensurate case in \cite{AlexandrovaMondie2021}.

%%%%%%%%%%%%%%%%%%%%%%%%%%%%%%%%%%%%%%%%%%%%%%%%%%%%%%%%%
\section{Ongoing and future research} \label{sec:OngoingandFuture}

Next, we outline current ongoing work and topics we believe are promising directions of research in the present framework.

%%%%%%%%%%%%%%%%%%%%%%%%%%%%%%%%%%%%%%%%%%%%%%%%%%%%%%%%%%%%%%
\subsection{Completion of current partial results}

Current research includes completing the results obtained up to now. The ultimate goal is a stability criterion in a finite number of operations for each class of systems in Section \ref{sec:state_of_art}. 
Generalizations or combinations of the cases described in Section \ref{sec:state_of_art} are possible. For example, neutral type systems with multiple delays in the difference operator and neutral type systems with distributed delays for which the construction of the Lyapunov matrix is available \citep{Aliseyko2019}.

Some generalizations deserve attention because they describe problems of practical interest. A first example is the combination of difference equation in continuous time with integral delay equation of the form
\begin{equation*} \label{eq:sys1}
	x(t)=\sum_{j=1}^m A_jx(t-h_j) +\int_{-H}^{0}G(\theta)x(t+\theta)\,d\theta,
\end{equation*}
which is an appropriate model for epidemics and population models. A second example is the interconnection of systems of differential delay equations with integral equations, which describes the closed-loop of state and input delay systems with predictor-based control law assigning a finite spectrum \citep{Kharitonov2014}. Here, a challenging problem is the distinct type of stability criteria for the differential and  integral equations.  

%%%%%%%%%%%%%%%%%%%%%%%%%%%%%%%%%%%%%%%%%%%%%%%%%%%%%%%%%
\subsection{Extensions to other classes of systems} \label{roadmap}

Based on the experience gained in the study of the cases in Section \ref{sec:state_of_art}, one can sketch a road map for presenting stability criteria in terms of the Lyapunov matrix of linear systems in general.

The prerequisite to be developed are:
\begin{itemize}
    \item determine the fundamental matrix and Cauchy formula for the system under consideration (provided that solutions of the system exist and are unique),\\
    \item define the concept of stability for the system,\\
    \item find, assuming the system stability, the form of the functional with prescribed quadratic negative derivative,\\
    \item define the corresponding Lyapunov matrix (notice that such choice is not unique),\\
    \item deduce properties of the Lyapunov matrix allowing its construction and develop the corresponding algorithm,\\
    \item present necessary and sufficient stability theorems in the form of a positivity condition for the obtained functional.
\end{itemize}

Once these elements are at hand, the steps presented in this contribution can be followed.
% \begin{itemize}
%     \item The initial function depending on the fundamental matrix that reveals a collection of necessary stability conditions in terms of the Lyapunov matrix must be determined.
%     \item The generalized properties allowing the reduction should be proved, without stability assumption.
%     \item For sufficiency, the error of approximation of arbitrary initial functions by the chosen initial function depending on the fundamental matrix must be estimated. 
%     \item The computation of all constants in the approximation error and in the stability theorem allow to charracterize the order of approximation for which the conditions are necessary and sufficient.
% \end{itemize}

We believe that the above prerequisites are available, or can be developed, in the cases of stochastic linear delay system with Wiener process, parabolic partial differential equations such as the heat equation, and hyperbolic ones such as the wave equation.

However, caution is in order: extensions are not straightforward as each class of systems presents challenges related to the Lyapunov matrix definition, its construction, and the stability theorems. 
% Indeed, the NP-hard nature of delay systems resurfaces, in particular in the approximation error for the initial function.
\subsection{Other approximations of the initial function}
Another strategy for presenting finite criterion of stability is possible if the prerequisites described in Section \ref{roadmap} are satisfied. Instead of using the initial functions in terms of the system fundamental matrix that reveal the Lyapunov matrix in the stability criterion, one can use an initial function expressed in terms of a given basis with free parameters. A bound on the  approximation error of arbitrary initial conditions of the compact set $\S$ with elements on this basis must be determined.

In this line of research, the initial function has been approximated by piecewise linear continuous functions in the case of systems with multiple delays in \cite{MedveZhabko2015}, and neutral type systems in \cite{ALEXANDROVA201983}. Legendre polynomials approximations were used for single-delay systems in \cite{bajodek:hal-03435028}. Other approximation choices of the initial functions are possible, for example, the Chebyshev polynomials used in \cite{Jarlebringetal2011} for the approximation of matrix $U$.

The order of approximation for which the necessary and sufficient stability condition is reached with these methods are substantially small compared to the approximation by fundamental matrices, thus making the results quite appealing. The reason is that the initial conditions $\psi_r$ defined in~\eqref{initial function definition}  poorly approximate arbitrary initial functions $\widehat{\ph}\in\Cf$ resulting in a large bound for the estimation error.

%However, within these approaches, the positivity test is carried out on the functional, demanding the computation of the integrals.

However, within the above mentioned approaches, the positivity test requires the computation of integrals. The proposed recursive formulae are time-consuming and prone to accumulating numerical errors. Regarding time consumption, there seems to be a trade-off between the quality of the approximations and the complexity of the recursive formula for computing the integrals. Notice also that results based on other approximations lose the simplicity of the criterion $\K_r>0$, which only requires computing the Lyapunov matrix at a finite number of points. 

%Notice also that results based on other approximations lose the simple expression of the criterion in terms of the delay Lyapunov matrix.

% \textcolor{red}{While it is true, I would not include the following paragraph. It is up to you.}

% \textcolor{red}{As a final observation, the combination of approaches would enhance the efficacy when testing a given parameter space: run first the low computational effort Lyapunov matrix-based test~\eqref{K_def} to discard unstable points with small values of $r$. For the remaining points, compute the order of piecewise or Legendre polynomials approximation ensuring sufficiency, and carry out the test.}

% \subsection{An open problem}
% Beyond the sleek dependence on the system Lyapunov matrix, a noteworthy feature of the stability condition $\K_r>0$ is that, according to what numerical examples show, it makes possible to rapidly approach the stability region in the space of parameters as $r$ increases. While the approach based on the compact set $\mathcal{S}$ and its approximation by a particular class of functions has allowed providing a number $\hat r$ for which we can ensure that the stability condition is indeed a stability criterion, the question of whether it is possible to find the exact number $r$ is still open. We believe that the approach on which the current results are grounded might not be appropriate for addressing such a question, but we do not discard the possibility of answering it.

%%%%%%%%%%%%%%%%%%%%%%%%%%%%%%%%%%%%%%%%%%%%%%%%%%%%%%%%%
\subsection{Alternative method for determining conditions in terms of the Lyapunov matrix}

The form of the necessary stability conditions in terms of the delay Lyapunov matrix can be uncovered directly. However, this approach relies on the stability assumption of the system; thus it cannot be used to prove the stability criterion. This strategy have been helpful in \cite{Rochaetal2016,delValleetal2018}. 

As the fundamental matrix $K(\cdot+\tau)$ is solution of~\eqref{eqSyst} for any $\tau\>0$, it satisfies the Cauchy formula.

\begin{lem} \label{Cauchy fundamental}
For $t\>0$ and $\tau\>0$, the fundamental matrix is such that
\begin{equation*}
	K(t+\tau) =K(t)K(\tau)+\sum_{j=1}^m\int_{-h_j}^{0} K(t-\theta-h_j)A_{1}K(\theta+\tau)\,d\theta.
\end{equation*}
\end{lem}

Consider~\eqref{initial function definition} as initial function. From equation~\eqref{Cauchy formula},
\begin{equation*}
	\begin{aligned}
		x(t,\psi_r) &=\sum_{i=1}^r\left(\vphantom{\int_{-h_j}^{0}} K(t)K(\tau_{i})\right. \\
        &\left. +\sum_{j=1}^m \int_{-h_j}^{0} K(t-\theta-h)A_{1}K(\theta +\tau_{i})\,d\theta\right)\gamma_{i}.
	\end{aligned}
\end{equation*}
In view of Lemma \ref{Cauchy fundamental}, this is
\begin{equation} \label{a}
	x(t,\psi_r)=\sum_{i=1}^{r}K(t+\tau_{i})\gamma_{i},\;\;t\> -H.
\end{equation}
The assumption of exponential stability of system~\eqref{eqSyst} and the substitution of~\eqref{a} into~\eqref{v1 integral form} yield
\begin{equation*}
	\begin{aligned}
	    &v_{1}(\psi_r) =\int_{-H}^{\infty} x^T(t,\psi_r)Wx(t,\psi_r)\,dt \\
	    &=\int_{-H}^{\infty}\left(\sum_{i=1}^{r}\gamma_{i}^T K^{T}(t+\tau_{i})\right) W\left(\sum_{j=1}^r K(t+\tau_j)\gamma_j\right)\,dt \\
	    &=\sum_{i=1}^{r}\sum_{j=1}^{r}\gamma_{i}^T \int_{-H}^{\infty} K^{T}(t+\tau_{i})WK(t+\tau_{j})\,dt\,\gamma_{j},
	\end{aligned}
\end{equation*}
and the change of variable $t+\tau_{i}=s$ gives
\begin{equation*}
	\begin{aligned}
		&v_{1}(\psi_r) =\sum_{i=1}^{r}\sum_{j=1}^{r}\gamma_{i}^T
        \left(\int_{\tau_{i}-H}^{\infty} K^T(s) WK(s+\tau_{j}-\tau_{i})\,ds\right)\gamma_{j} \\
		&=\sum_{i=1}^{r}\sum_{j=1}^{r}\gamma_{i}^T \left(\int_{0}^{\infty} K^{T}(s)W K(s+\tau_{j}-\tau_{i})\,ds\right) \gamma_{j}.
	\end{aligned}
\end{equation*}
By using equality~\eqref{ec:defiU}, we get
\begin{equation*}
	v_{1}(\psi_r)=\sum_{j=1}^{r}\sum_{i=1}^{r}\gamma_{i}^T U(\tau_{j}-\tau_{i})\gamma_{j} =\gamma ^{T}\left[\vph U(\tau_{j}-\tau_{i})\right]_{i,j=1}^r \gamma,
\end{equation*}
where $\gamma=\left(\gamma_1^T\,\ldots\,\gamma_r^T\right)^T$. From Theorem \ref{KharZabMod} we conclude that
\begin{equation*}
	\left[U(\tau_j-\tau_i)\vph \right]_{i,j=1}^r > 0.
\end{equation*}

Current research aims at the extension of this result to unstable systems. The idea is to consider parametrized form of functional $v_1$ and of the Lyapunov matrix $U$:
\begin{gather*}
    \widehat v_1(\ph,p) =\int_{-H}^{\infty} e^{-pt} x^T(t,\ph)Wx(t,\ph)\,dt, \\
    \widehat U(\tau,p) =\int_0^{\infty} e^{-pt} K^T(t)WK(t+\tau)\,dt.
\end{gather*}
Parameter $p$ is complex. It is clear that such integrals converge only for $p$ with sufficiently large real part. But for other values of $p$ we can use the analytic continuation. Thus, we can prove all the desired properties for the improper integrals and then continue them analytically.

One can show that the analytic continuation to point $p=0$ exists, if the Lyapunov condition holds.

%%%%%%%%%%%%%%%%%%%%%%%%%%%%%%%%%%%%%%%%%%%%%%%%%%%%%%%%%%%%%%%%%
\section{Conclusion}

We have presented sleek necessary and sufficient stability conditions that match the delay-free case: they are expressed as a positivity condition of the system Lyapunov matrix and can be verified in a finite number of operations. The general framework we introduce allows the presentation of other finite criteria of stability and extensions to other classes of linear systems.

We hope the detailed presentation of the multiple delay case corresponding to papers scattered in the past decade's literature will help to understand the topic and encourage future research in this direction.

\appendix

\section{Proof of necessary stability conditions in terms of functional $v_1$}

\subsection{Proof of Theorem \ref{KharZabMod}} \label{app:KharZabMod}

Consider the functional
\begin{equation*}
    \begin{aligned}
        \widetilde{v}(\ph)&=v_1(\ph)-\frac{1}{m+1}\sum_{j=0}^m \int_{-H}^{-h_j}\ph^T(\theta)W\ph(\theta)\,d\theta \\
        &-\alpha_0\|\ph(0)\|^2.
    \end{aligned}
\end{equation*}
Its derivative along the solutions of system~\eqref{eqSyst} is
\begin{equation*}
    \begin{aligned}
        \D\widetilde v(\ph)&=\D v_1(\ph) -\frac{1}{m+1}\sum_{j=0}^m \ph^T(-h_j)W \ph(-h_j) \\
        &+\ph^T(-H)W \ph(-H) -2\alpha_0 \ph^T(0) \sum_{j=0}^mA_j\ph(-h_j).
    \end{aligned}
\end{equation*}
By~\eqref{eqDv1},
\begin{equation} \label{eqNecCondFunctAux}
    \begin{aligned}
        &\D\widetilde v(\ph)= -\frac{1}{m+1}\sum_{j=0}^m \ph^T(-h_j)W \ph(-h_j) \\
        &-2\alpha_0\ph^T(0)\sum_{j=0}^mA_j\ph(-h_j).
    \end{aligned}
\end{equation}
As the first term is negative definite, there exists sufficiently small positive $\alpha_0$ to preserve negativity. For such value of $\alpha_0$
\begin{equation*}
    \D \widetilde v\left(\ph\right)\<0.
\end{equation*}
As the systems is exponentially stable, i.\,e., $\lim\limits_{t\to\infty} \widetilde v(x_t(\ph)) =0$, then
\begin{equation*}
    \widetilde v(\ph) =-\int_0^{\infty} \D \widetilde v\left(x_t(\ph)\right)\,dt \>0.
\end{equation*}
Therefore,
\begin{equation*}
    \begin{aligned}
        v_1(\ph)& \>\frac{1}{m+1}\sum_{j=0}^m \int_{-H}^{-h_j}\ph^T(\theta)W\ph(\theta)\,d\theta +\alpha_0\|\ph(0)\|^2 \\
        &\>\frac{1}{m+1}\int_{-H}^0\ph^T(\theta)W\ph(\theta)\,d\theta +\alpha_0\|\ph(0)\|^2.
    \end{aligned}
\end{equation*}
We can take $\alpha_1=\dfrac{\lambda_{\min}(W)}{m+1}$.

%%%%%%%%%%%%%%%%%%%%%%%%%%%%%%%%%%%%%%%%%%%%%%%%%%%%%%%%%%%%%%%%%
\subsection{Proof of Theorem \ref{theo:stability_v1_set}} \label{app:theo:stability_v1_set}

An obvious corollary of the proof of Theorem \ref{KharZabMod}:
\begin{equation*}
    v_1(\ph) \>\alpha_0,\;\;\ph\in\S.
\end{equation*}
Let us show that $\D\widetilde v$ in the proof of Theorem \ref{KharZabMod} is non-positive for $\alpha_0=\alpha_0^{\star}$.

Equality~\eqref{eqNecCondFunctAux} can be rewritten as
\begin{equation*}
    \begin{aligned}
        \D\widetilde v(\ph)&= -\frac{1}{m+1}\left(\sum_{j=0}^m \ph^T(-h_j)W \ph(-h_j)\right. \\
        &\left.-2\dfrac{1}{\lambda_{\min}(P)}\ph^T(0) \sum_{j=0}^mA_j\ph(-h_j)\right) \\
        &= -\frac{1}{m+1}\left(\vphantom{\dfrac{1}{\lambda_{\min}(P)}} \widetilde\ph^T \cdot I_{m+1}\otimes W \cdot\widetilde\ph\right. \\
        &\left.-2\dfrac{1}{\lambda_{\min}(P)}\widetilde\ph^T\cdot e_1\otimes I_n \cdot A \cdot \widetilde\ph\right),
    \end{aligned}
\end{equation*}
where $\widetilde\ph=\left(\ph^T(0),\ph^T(-h_1),\ldots,\ph^T(-h_m)\right)^T$.

Symmetrization leads to
\begin{equation} \label{eqNecCondFunctAux2}
    \begin{aligned}
        \D\widetilde v(\ph)&= -\frac{1}{(m+1)\lambda_{\min}(P)}\widetilde\ph^T \left(\vph\lambda_{\min}(P) I_{m+1}\otimes W \right. \\
        &\left.\vph - e_1\otimes I_n\cdot A -A^T\cdot e_1^T\otimes I_n\right)\widetilde\ph.
    \end{aligned}
\end{equation}
As $W$ is positive definite, we can apply Cholesky decomposition:
\begin{equation*}
    I_{m+1}\otimes W = I_{m+1}\otimes C^T \cdot I_{m+1}\otimes C,
\end{equation*}
where $C$ is invertible. By Schur complement, the eigenvalues of $P$ are the same as those of matrix
\begin{equation*}
    \widetilde P=I_{m+1}\otimes C^{-T} \cdot \left(e_1\otimes I_n\cdot A + A^T\cdot e_1^T\otimes I_n\right) \cdot I_{m+1}\otimes C^{-1}.
\end{equation*}

Now~\eqref{eqNecCondFunctAux2} takes the form
\begin{equation*}
    \begin{aligned}
        \D\widetilde v(\ph)&= -\frac{1}{(m+1)\lambda_{\min}(P)}\widetilde\ph^T\cdot I_{m+1}\otimes C^T \\
        &\cdot\left(\vph\lambda_{\min}(P) I_{n(m+1)}-\widetilde P\right)I_{m+1}\otimes C\cdot\widetilde\ph.
    \end{aligned}
\end{equation*}
The matrix in brackets is negative semidefinite. It remains two show that $\lambda_{\min}(P)<0$.

If, on the contrary, all eigenvalues of $P$ (equivalently, $\widetilde P$) are nonnegative, then
\begin{equation} \label{eqPneg}
    \mu^T\left(e_1\otimes I_n\cdot A + A^T\cdot e_1^T\otimes I_n\right)\mu\>0
\end{equation}
for any $\mu\in\R^{n(m+1)}$.

Assume for simplicity that matrix $A_m$ is nontrivial. Let $\gamma\in\R^n$ be such that $A_m\gamma\neq 0$. Take
\begin{equation*}
    \mu=-e_1\otimes \left(\delta A_m\gamma\right) +e_{m+1}\otimes \gamma,
\end{equation*}
where $e_{m+1}=(0,\ldots,0,1)^T\in\R^{m+1}$, number $\delta>0$. Substitution of the above value of $\mu$ into~\eqref{eqPneg} and simple computations give
\begin{equation*}
    -2\delta \gamma^T A_m^TA_m\gamma +\delta^2 \gamma^TA_m^T\left(A_0^T+A_0\right)A_m\gamma\>0.
\end{equation*}
For sufficiently small $\delta$, this inequality fails. This contradiction finishes the proof.

%%%%%%%%%%%%%%%%%%%%%%%%%%%%%%%%%%%%%%%%%%%%%%%%%%%%%%%%%%%%%%%%%
\section{Proof of necessary instability conditions in terms of functional $v_1$}

\subsection{Proof of Theorem \ref{ZhMedvMod}} \label{app:ZhMedvMod}

As system~\eqref{eqSyst} is unstable and the Lyapunov condition holds, there exists an eigenvalue $s_0=p+iq$ with $p>0$ and $q\>0$, and two vectors $C_1,\:C_2\in\R^n$ such that 
\begin{equation*} \label{ec:IL_sol_unstable}
	x(t,\widehat\ph)=\ell e^{p t}\phi(t),\;\;t\in(-\infty,\infty),
\end{equation*}
where $\ell$ is an arbitrary real number and $\phi(t)=\cos\left(qt\right) C_1-\sin\left(qt\right) C_2$, is solution of~\eqref{eqSyst} corresponding to the initial function
\begin{equation*}
	\widehat\ph(\theta)=x(\theta,\widehat\ph),\:\theta\in [-H,0].
\end{equation*}
Moreover, the vectors $C_1$ and $C_2$ satisfy $C_1^TC_2=0$, $\|C_1\|=1$ and $\|C_2\|\< 1$, cf. \cite[Lemma 8]{Gomezetal2019}.

Since for any $\ph\in\PC$
\begin{equation*}
	\D v_1(\ph)=-\ph^T(-H)W\ph(-H),
\end{equation*}
we have that
\begin{equation*}
    \begin{aligned}
	    v_1(x_T(\widehat\ph)) -v_1(\widehat\ph)&=-\int_0^T x^T(t-H,\widehat\ph) Wx(t-H,\widehat\ph)\,dt \\
	    &=-\int_{-H}^{T-H} x^T(t,\widehat\ph) Wx(t,\widehat\ph)\,dt,
	\end{aligned}
\end{equation*}
where $W$ is positive definite. If $T=2\pi/q$ for $q\neq 0$, and $T=1$ for $q=0$, it is a period of function $\phi$. Therefore, we have $x_T(\widehat\ph)=e^{p T}\widehat\ph$ and
\begin{equation*}
	v_1(x_T(\widehat\ph))=e^{2p T}v_1(\widehat\ph),
\end{equation*}
which implies that
\begin{equation*}
    \begin{aligned}
	    v_1(\widehat\ph)&=-\dfrac{1}{e^{2p T}-1}\int_{-H}^{T-H} x^T(t,\widehat\ph)W x(t,\widehat\ph)\,dt \\
	    &\< -\dfrac{\lambda_{\min}(W)}{e^{2p T}-1} \int_{-H}^{T-H}\|x(t,\widehat\ph)\|^2\,dt.
	\end{aligned}
\end{equation*}
By the above mentioned properties of vectors $C_1$ and $C_2$,
\begin{equation*}
    \begin{aligned}
        \int_{-H}^{T-H}\|x(t,\widehat\ph)\|^2\,dt &\>\ell^2\int_{-H}^{T-H}e^{2p t}\cos^2(q t)\,dt \\
        &=\ell^2 \dfrac{e^{-2pH}\(e^{2pT}-1\)}{4p} f(q),
	\end{aligned}
\end{equation*}
where
\begin{equation*}
	f(q)=\cos^2(qH)+\dfrac{\(\vphS p\cos(qH) -q\sin(qH)\)^2}{p^2+q^2}.
\end{equation*}
Hence,
\begin{equation*}
	v_1(\widehat\ph)\< -\ell^2 \dfrac{\lambda_{\min}(W)e^{-2pH}}{4p} f(q).
\end{equation*}

It is easy to see that $f(q)>0$ for any $q$. Therefore, for any $\beta>0$ there exists $\ell$ such that $v_1(\widehat\ph)\<-\beta$.

%%%%%%%%%%%%%%%%%%%%%%%%%%%%%%%%%%%%%%%%%%%%%%%%%%%%%%%%%%%%%%
\subsection{Proof of Theorem \ref{theo:instability}} \label{app:instability_theorem}

The proof is based on the one of Theorem \ref{ZhMedvMod}. We take the same solution $x(\cdot,\widehat\ph)$ with $\ell=1$. For the corresponding $\widehat\ph$
\begin{equation*}
	v_1(\widehat\ph)\< -\dfrac{\lambda_{\min}(W)e^{-2p H}}{4p} f(q).
\end{equation*}

We focus now on providing a strictly positive lower bound for the function $f$. On the one hand
\begin{equation*}
	f(q)\> \cos^2(qH).
\end{equation*}
On the other hand
\begin{equation*}
    \begin{aligned}
	    f(q)&=1+\dfrac{p}{\sqrt{p^2+q^2}} \\
	    &\cdot\left(\dfrac{p}{\sqrt{p^2+q^2}}\cos(2qH) -\dfrac{q}{\sqrt{p^2+q^2}} \sin(2qH)\right)\\
	    &\> 1-\dfrac{p}{\sqrt{p^2+q^2}}.
    \end{aligned}
\end{equation*}
As $a$ is an upper estimate of the real part of the
rightmost root of system (\ref{eqSyst}), $p\<a$ and
\begin{equation*}
    f(q)\>1-\dfrac{a}{\sqrt{a^2+q^2}}.
\end{equation*}

Thus, to find a positive lower bound for the function $f$ independent of $q$, one needs to solve the following optimization problem:
\begin{equation*}
	\min_{q\in[0,\infty)} \max\left\{\cos^2(qH),\,1-\dfrac{a}{\sqrt{a^2+q^2}} \right\}.
\end{equation*}
It is easy to see that the minimum is achieved at point $q_0\in\left(0,\dfrac{\pi}{2H}\right)$, such that
\begin{equation*}
	\cos^2(q_0H)=1-\dfrac{a}{\sqrt{a^2+q_0^2}}.
\end{equation*}
The substitution $b=q_0H$ and some simple transformations lead to equality~\eqref{ec:b}. Since
\begin{equation*}
	g(b)=((aH)^2+b^2)\sin^4(b)-(aH)^2, \;\;b\in\left[0,\dfrac{\pi}{2}\right],
\end{equation*}
is an increasing function,
\begin{equation*}
	g(0)=-(aH)^2<0\:\:\:\text{and}\:\:\:g \left(\frac{\pi}{2}\right) =\dfrac{\pi^2}{4}>0,
\end{equation*}
the number $b$ on $(0,\pi/2)$ such that $g(b)=0$ exists and is unique. Therefore, by the fact that $0<p\<a$, we get the desired result:
\begin{equation*}
	v_1(\widehat\ph) \<-\dfrac{\lambda_{\min}(W)e^{-2aH}}{4a}\cos^2(b).
\end{equation*}

It remains to show that $\widehat\ph\in\S$. By construction, for $t\in(-\infty,\infty)$
\begin{equation*}
    \|x(t,\widehat\ph)\|^2=e^{2pt}\left(\cos^2(qt) +\sin^2(qt)\|C_2\|^2\right).
\end{equation*}
On $(-\infty,0]$ this value is not greater than $1$ and achieves $1$ at zero. In particular, $\|\widehat\ph(0)\|=1$, $\|\widehat\ph(\theta)\|\<1$, $\theta\in(-H,0)$.

Also, for $t\in(-\infty,0)$
\begin{equation*}
    \begin{aligned}
        \|\dot x(t,\widehat\ph)\| &=\left\|\sum_{j=0}^m A_j x(t-h_j,\widehat\ph)\right\| \\
        &\<\sum_{j=0}^m \|A_j\| \|x(t-h_j,\widehat\ph)\|\<\sum_{j=0}^m \|A_j\|=M.
    \end{aligned}
\end{equation*}
Hence, $\|\widehat\ph'(\theta)\|\<M$, $\theta\in[-H,0]$.

%%%%%%%%%%%%%%%%%%%%%%%%%%%%%%%%%%%%%%%%%%%%%%%%%%%%%%%%%%%%%%%%%
\section{Proof of Lemma \ref{lemma Z reduction}} \label{app:lemma Z reduction}

To reduce the expression of the functional valued at the special initial function to expression~\eqref{z functional}, the following technical property, which holds true for $\tau_1\>0$, $\tau_2\in\R$ is needed (see, \cite{EgorovMondie2014} for the proof):
\begin{equation} \label{GAPint}
    \begin{aligned}
	    U(\tau_{1}&+\tau_{2})=U(\tau_{2})K(\tau_{1}) \\
	    &+\sum_{j=1}^{m}\int_{-h_{j}}^{0}U(\tau_{2}-\theta-h_{j}) A_{j}K(\tau_{1}+\theta)\,d\theta \\
	    &+\int_{-H}^{0}K^{T}(-\tau_{2}+\theta)WK(\tau_{1}+\theta)\,d\theta.
    \end{aligned}
\end{equation}
A direct consequence of this formula for $\tau_1\>0$, $\tau_2\>0$ is
\begin{equation} \label{CauchyLGen}
    \begin{aligned}
	    U(\tau_1&+\tau_2)=U(\tau_2)K(\tau_1) \\
	    &+\sum_{j=1}^m\int_{-h_j}^0U(\tau_2-\theta-h_j) A_j K(\tau_1+\theta)\,d\theta.
    \end{aligned}   
\end{equation}

The substitution of~\eqref{phipsi} into~\eqref{bilinear functional} yields
\begin{equation*}
    \begin{aligned}
        &z(K(\tau_1+\cdot)\mu,K(\tau_2+\cdot)\eta) =\mu^T \left[\vphantom{\int_{-H}^0}K^T(\tau_1)U(0)K(\tau_2)\right. \\
        &+K^T(\tau_1)\sum_{j=1}^m\int_{-h_j}^0 U(-\theta-h_j)A_jK(\tau_2+\theta)\,d\theta \\
        &+\sum_{k=1}^m\int_{-h_k}^0K^T(\tau_1+\theta)A_k^T U(\theta+h_k)\,d\theta K(\tau_2) \\
        &+\sum_{k=1}^m\int_{-h_k}^0K^T(\tau_1+\theta_1) A_k^T \sum_{j=1}^m\int_{-h_j}^0U(\theta_1+h_k-\theta_2-h_j)A_j \\
        &\left.\cdot K(\tau_2+\theta_2)\,d\theta_2\,d\theta_1 +\int_{-H}^0 K^T(\tau_1+\theta)WK(\tau_2+\theta) \,d\theta\right]\eta.
    \end{aligned}
\end{equation*}
Equality~\eqref{CauchyLGen}, applied to the first pair of terms and to the next one, implies that
\begin{equation*}
    \begin{aligned}
        &z(K(\tau_1+\cdot)\mu,K(\tau_2+\cdot)\eta)= \mu^T \left[\vphantom{\int_{-H}^0} K^T(\tau_1)U(\tau_2)\right. \\
        &+\sum_{k=1}^m\int_{-h_k}^0 K^T(\tau_1+\theta)A_k^T U(\tau_2+\theta+h_k)\,d\theta \\
        &+\left.\int_{-H}^0 K^T(\tau_1+\theta)W K(\tau_2+\theta)\,d\theta\right]\eta.
    \end{aligned}
\end{equation*}
Property~\eqref{GAPint} leads to the desired result.

\bibliography{autosam}           % and a bib file to produce the 
                                 % bibliography (preferred). The
                                 % correct style is generated by
                                 % Elsevier at the time of printing.

%\begin{thebibliography}{99}     % Otherwise use the  
                                 % thebibliography environment.
                                 % Insert the full references here.
                                 % See a recent issue of Automatica 
                                 % for the style.
%  \bibitem[Heritage, 1992]{Heritage:92}
%     (1992) {\it The American Heritage. 
%     Dictionary of the American Language.}
%     Houghton Mifflin Company.
%  \bibitem[Able, 1956]{Abl:56}
%     B.~C.~Able (1956). Nucleic acid content of macroscope. 
%     {\it Nature 2}, 7--9. 
%  \bibitem[Able {\em et al.}, 1954]{AbTaRu:54}   
%     B.~C. Able, R.~A. Tagg, and M.~Rush (1954).
%     Enzyme-catalyzed cellular transanimations.
%     In A.~F.~Round, editor, 
%     {\it Advances in Enzymology Vol. 2} (125--247). 
%     New York, Academic Press.
%  \bibitem[R.~Keohane, 1958]{Keo:58}
%     R.~Keohane (1958).
%     {\it Power and Interdependence: 
%     World Politics in Transition.}
%     Boston, Little, Brown \& Co.
%  \bibitem[Powers, 1985]{Pow:85}
%     T.~Powers (1985).
%     Is there a way out?
%     {\it Harpers, June 1985}, 35--47.

%\end{thebibliography}

\end{document}